\documentclass[11pt,reqno,a4paper]{amsart}
\usepackage{color}
\usepackage{amssymb,amsmath}
\usepackage[latin1]{inputenc}
\usepackage[active]{srcltx}
\usepackage{graphicx}
\usepackage{exscale}
\usepackage{textgreek}
\usepackage{epsfig,graphics}
\usepackage{psfrag}
\usepackage{caption}
\usepackage{subcaption}
\def\N{\mathbb{N}}
\def\R{\mathbb{R}}
\def\m1{{I\!\!M}}

\renewcommand{\to}{\rightarrow}
\newcommand{\pa}{\partial}

\newcommand{\ino}{\int_{\Omega}}

\newcommand{\ainf}{\mbox{as\;}\;n\to+\infty}

\newcommand{\fo}{\forall}

\newcommand{\bop}{\mbox{\rm o}}

\newcommand{\rife}[1]{(\ref{#1})}
\newcommand{\ov}[1]{\overline{#1}}
\newcommand{\un}[1]{\underline{#1}}

\newcommand{\sscp}{\scriptscriptstyle}
\newcommand{\dsp}{\displaystyle}
\newcommand{\noi}{\noindent}
\renewcommand{\dfrac}{\displaystyle\frac}
\newcommand{\finedim}{\hspace{\fill}$\square$}
\newcommand{\intbar}{\mathop{\int\makebox(-15.5,0){\rule[6pt]{.7em}{0.3pt}}%
\kern-6pt}\nolimits}

\newcommand{\ii}{\infty}

\newcommand{\eps}{\varepsilon}
\newcommand{\dt}{\delta}

\newcommand{\al}{\alpha}

\newcommand{\sg}{\sigma}

\newcommand{\om}{\Omega}
\newcommand{\lm}{\lambda}

\newcommand{\e}[1]{{\,\dsp e^{\dsp #1}}}

\newcommand{\rl}{\mbox{\Large \textrho}_{\!\sscp \lm}}

\renewcommand{\rho}{\mbox{\Large \textrho}}
\newcommand{\rh}{\mbox{\Large \textrho}}

\newcommand{\pl}{\psi_{\sscp \lm}}

\newcommand{\ple}{\psi_{\sscp \le}}
\newcommand{\vb}{v_{\sscp \lm}}

\newcommand{\ssb}{\sscp \lm}
\newcommand{\sse}{\sscp \mathcal{E}}

\renewcommand{\le}{\lm(E)}

\newcommand{\ol}{\ssb,\sscp 0}

\newcommand{\prl}{{\textbf{(}\mathbf P_{\mathbf \lm}\textbf{)}}}

\newtheorem{theorem}{Theorem}[section]
\newtheorem{proposition}[theorem]{Proposition}
\newtheorem{lemma}[theorem]{Lemma}
\newtheorem{corollary}[theorem]{Corollary}
\newtheorem{remark}[theorem]{Remark}
\newtheorem{definition}[theorem]{Definition}
\newcommand{\brm}{\begin{remark}\rm}
\newcommand{\erm}{\end{remark}}
\newcommand{\bdf}{\begin{definition}\rm}
\newcommand{\edf}{\end{definition}}
\newcommand{\bte}{\begin{theorem}}
\newcommand{\ete}{\end{theorem}}
\newcommand{\bpr}{\begin{proposition}}
\newcommand{\epr}{\end{proposition}}
\newcommand{\ble}{\begin{lemma}}
\newcommand{\ele}{\end{lemma}}
\newcommand{\bco}{\begin{corollary}}
\newcommand{\eco}{\end{corollary}}
\newcommand{\beq}{\begin{equation}}
\newcommand{\eeq}{\end{equation}}
\newcommand{\bdm}{\begin{displaymath}}
\newcommand{\edm}{\end{displaymath}}

\newcommand{\graf}[1]{\left\{\begin{array}{ll}#1\end{array}\right.}

\def\sideremark#1{\ifvmode\leavevmode\fi\vadjust{\vbox to0pt{\vss
 \hbox to 0pt{\hskip\hsize\hskip1em \vbox{\hsize2.1cm\tiny\raggedright\pretolerance10000 \noindent #1\hfill}\hss}\vbox to15pt{\vfil}\vss}}}

\begin{document}
\numberwithin{equation}{section}
\parindent=0pt
\hfuzz=2pt
\frenchspacing

\title[Mean field equations and Onsager's theory]{Global bifurcation analysis of mean field equations and
the Onsager microcanonical description of two-dimensional turbulence}

\author[D.B.]{Daniele Bartolucci$^{(1,\dag)}$}


\thanks{$^{(1)}$Daniele Bartolucci, Department of Mathematics, University
of Rome {\it "Tor Vergata"}, \\  Via della ricerca scientifica n.1, 00133 Roma,
Italy. e-mail:bartoluc@mat.uniroma2.it}

\thanks{$^{(\dag)}$Research partially supported by FIRB project "{\em
Analysis and Beyond}",  by PRIN project 2012, ERC PE1\_11,
"{\em Variational and perturbative aspects in nonlinear differential problems}", and by the Consolidate the Foundations
project 2015 (sponsored by Univ. of Rome "Tor Vergata"),  ERC PE1\_11,
"{\em Nonlinear Differential Problems and their Applications}"}

\begin{abstract} On strictly starshaped domains of second kind we find natural sufficient conditions which allow the
solution of two long standing open problems closely related to the mean field equation $\prl$ below.\\
 On one side we catch the \un{global} behaviour of the Entropy
for the mean field Microcanonical Variational Principle ((MVP) for short)
arising in the Onsager description of two-dimensional turbulence.
This is the completion of well known results first established in \cite{clmp2}.
Among other things we find a full unbounded interval of strict convexity of the Entropy.\\
On the other side, to achieve this goal,  we have to provide a detailed qualitative description of
the \un{global} branch of solutions of $\prl$ emanating from $\lm=0$ and crossing
$\lm=8\pi$. This is the completion of well known results first established in \cite{suz} and \cite{CCL} for $\lm\leq 8\pi$, and it has 
an independent mathematical interest, since the shape of global branches of semilinear elliptic equations, with very few well known
exceptions, are poorly understood.
It turns out that the (MVP) suggests the right variable (which is the energy) to be used to
obtain a global parametrization of solutions of $\prl$. A crucial spectral simplification
is obtained by using the fact that, by definition, solutions of the (MVP) maximize the entropy at fixed energy and total vorticity.
\end{abstract}
\maketitle

\tableofcontents

\section{\bf Introduction.}
\setcounter{equation}{0}
We are concerned with two long standing open problems closely related to the mean field equation $\prl$ below on
smooth, open, bounded and connected domains $\om\subset\R^2$. For a large class of domains
(known as domains of second kind \cite{clmp2}, see Definition \ref{def1}), and with the exception of the results in \cite{clmp2},
we miss the details of the mean field supercritical thermodynamics of the vorticity distribution (whose equilibrium potentials satisfy $\prl$ with $\lm>8\pi$)
in the Onsager's statistical description of two dimensional turbulence \cite{On}.
On the other side, still for domains of second kind, and with the exception of the results in \cite{BLin3}, \cite{CCL} and more recently in
\cite{BdM2}, we don't know much about the possible continuation of the branch of solutions of $\prl$ emanating from $\lm=0$ and crossing $\lm=8\pi$.
The lack of our understanding of the latter is part of the cause of the lack of our understanding of
the former problem, and  it is not surprising that both problems are essentially solved at once, as we will
discuss below. As a matter of fact, the physical problem suggests the right variable, which is the energy, to be used to
obtain a \un{global} parametrization of solutions of $\prl$.

\bigskip

Let us define,
$$
\mathcal{P}_{\sscp \om}=\left\{\rho\in L^{1}(\om)\,|\,\rho\geq 0\;\mbox{a.e. in}\;\om,\;\ino\rho =1 \right\},
$$
and let $G(x,y)$ be the unique  solution of,
\beq\label{Green}
\graf{
-\Delta G(x,y)= \delta_{y}&
\mbox{in}\hspace{.4cm}\;\; \om, \\
\hspace{0.7cm}G(x,y)=0 &\mbox{on}\hspace{.1cm}\;\; \pa \om.}
\eeq

We define the entropy and the energy of  the vorticity density $\rho\in \mathcal{P}_{\sscp \om}$ as,

$$
\mathfrak{S}(\rho)=-\ino \rh\log(\rh),\quad
\mathcal{E}(\rho)=\frac12 \ino \rho G[\rho],
$$
respectively, where,
$$
G[\rho](x)=\ino G(x,y)\rho(y)\,dy,\;x\in\om.
$$

For any  $E\in(0,+\infty)$ we consider the Microcanonical Variational Principle,
$$
S(E)=\sup \left\{ \mathfrak{S}(\rho),\; \rho \in \mathcal{M}_E \right\},\quad
\mathcal{M}_E=\left\{\rho\in\mathcal{P}_{\sscp \om}\,|\,\mathcal{E}(\rho)=E\right\}.\qquad \mbox{\bf (MVP)}
$$

For $\psi\in H^{1}_0(\om)$ and $\lm\in\R $ let us define,
$$
\rl(\psi) =\dfrac{\dsp \e{\lm \psi}}{\dsp \ino \e{\lm \psi}}.
$$

By the Moser-Trudinger \cite{moser} and Jensen's inequalities we see that $\rl(\psi)\in \mathcal{P}_{\sscp \om}$. It has been shown in \cite{clmp2} that
the densities solving the ${\bf (MVP)}$, take the form $\rh=\rl(\psi)$ and satisfy, for some $\lm\in \R$, the Euler-Lagrange equations,
$$
\graf{-\Delta \psi =\dfrac{\dsp \e{\lm \psi}}{\dsp \ino \e{\lm \psi}}\quad \om \\ \psi = 0 \qquad \qquad\pa\om}\qquad\qquad \qquad \qquad \prl
$$
\brm{\it
More exactly $\psi$ {\rm (}the stream function{\rm)} satisfies {\rm$\psi=G[\rl(\psi)]$},
which implies, by the results in \cite{bm}, that $\psi\in L^{\ii}(\om)$.
Therefore, since $\om$ is smooth, then standard elliptic regularity theory shows that $\psi$
is of class $C^{2,\alpha}_0(\om)\cap C^1(\ov{\om})$.
In particular, we conclude that {\rm $\rl(\psi)$} is uniformly bounded from above and from below away from zero.
We will often use this properties in the rest of this paper without further comments.}
\erm

Let us recall that, by the Jensen's inequality, we have, $\mathfrak{S}(\rl)\leq\log(|\om|),$ where the equality holds if and only if
$\rl(x)\equiv \rh_0(x)=\frac{1}{|\om|},\,\fo x\in\ov{\om}$. Therefore $\rh_0$ is the unique solution of the ${\bf (MVP)}$ at energy $E=E_0$, that is,
$$
E_0:=\mathcal{E}\left(\rl(\pl)\right)\left.\right|_{\lm=0}=\mathcal{E}\left(\rh_0\right)=\frac{1}{|\om|^2}\ino \ino G(x,y) dx dy.
$$

Among many other things which we will not discuss here, the following facts has been proved in \cite{clmp2}:\\

\textbf{(MVP1)}: for each $E> 0$ there exists $\le\in\R$ and a solution $\pl=\ple$ of $\prl$ such that $\rho_{\sscp \le}(\ple)$
solves the ${\bf (MVP)}$. Actually $\prl$ is the Euler-Lagrange variational equation associated to the ${\bf (MVP)}$,
the energy constraint being imposed by choosing $\lm=\le$ such that $\mathcal{E}\left(\rl(\pl)\right)=E$.
The constraint of unitary mass for the density
is already imposed by the expression of $\rl(\psi)$.\\

\textbf{(MVP2)}: $S(E)$ is continuous and
\beq\label{conc}
\rh_{\sscp \lm(E)}(\ple)\rightharpoonup \dt_{x_0},\quad \mbox{as}\quad E\to +\infty,
\eeq
where $x_0$ is a maximum point of $\gamma(x):=\gamma(x,x)$, and $\gamma(x,y)=G(x,y)+\frac{1}{2\pi}\log(|x-y|)$.\\

\bdf\label{def1} {\it
Let $\om\subset\R^2$ be any smooth and bounded domain. We say that $\om$ is of {\bf second kind} if {\rm $\prl$} admits a solution for $\lm=8\pi$. Otherwise
$\om$ is said to be of {\bf first kind}.}
\edf

The definition of domains of first/second kind was first introduced in \cite{clmp1} with an equivalent but different formulation.
We refer to \cite{CCL} and \cite{BLin3} for a complete
discussion about the characterization of domains of first/second kind and related examples.
Among many other things, it has been proved in \cite{BdM2} that
there exists a universal constant $I>4\pi$ such that any convex domain whose isoperimetric ratio $I(\om)$ satisfies $I(\om)\geq I$ is of second kind.\\
It is well known that solutions of $\prl$ exists and are unique for any $\lm<8\pi$ and are unique whenever they exist for $\lm=8\pi$,
see either \cite{suz} for $\lm<8\pi$ and $\om$ simply connected or \cite{CCL} for $\lm\leq 8\pi$ and $\om$ simply connected
and more recently \cite{BLin3} for $\lm\leq 8\pi$ and $\om$ any connected domain.
\bdf\label{def2} {\it
If $\om$ is of {\bf second kind} we define  $E_{8\pi}:=\mathcal{E}\left(\mbox{\rm \rh}_{\sscp 8\pi}(\psi_{\sscp 8\pi})\right)$.
Otherwise, if $\om$ is of {\bf first kind}, we set $E_{8\pi}=+\infty$.}
\edf

The mean field thermodynamics of the system is well understood for $E<E_{8\pi}$ in terms of the unique solutions of $\prl$ for $\lm<8\pi$.
In this range of energies  it has been proved in \cite{clmp2} that:\\

\textbf{(MVP3)}: if $\om$ is simply connected, then for each $E<E_{8\pi}$ there exists $\lm=\lm(E)\in (-\infty,8\pi)$ such that $\lm(E)$ is well defined, continuous
and monotonic increasing, $\lm(E)\nearrow 8\pi$ as $E \nearrow  E_{8\pi}$ and $S(E)=\mathfrak{S}\left(\rho_{\sscp \le}(\ple)\right)$ is
smooth and concave in $(0,E_{8\pi})$. This result has been recently generalized to cover the case
where $\om$ is any bounded and connected domain in \cite{BLin3}.

\brm {\it From the physical point of view the natural variable arising as the Lagrange multiplier for the Energy constraint is not $\lm$ but $\beta=-\lm$,
which would correspond to the inverse statistical temperature $\beta_{stat}=\frac{1}{\kappa T_{stat}}$,
where  $\kappa$ is the Boltzmann constant. According to general principles of statistical thermodynamics, whenever the entropy $S$ happens to be
differentiable at the value of the energy $E$, we define,
$$
\beta(E):=\frac{dS(E)}{d E}.
$$
}
\erm

It has been shown in  \cite{CCL} and \cite{BLin3} that $\om$ is of second kind if and only if the branch of unique solutions of $\prl$
with $\lm<8\pi$ is uniformly
bounded, that is, $\om$ is of second kind if and only if $\pl$ converges smoothly to a solution of $\prl$ with $\lm=8\pi$, say $\psi_{\sscp 8\pi}$,
as $\lm\nearrow 8\pi$.  As a consequence,
in this situation it is not difficult to see that $S(E)$ extends by continuity from $(0,E_{8\pi})$ to $(0,E_{8\pi}]$ where
$S(E_{8\pi})=\mathfrak{S}\left(\rho_{\sscp 8\pi}(\psi_{\sscp 8\pi})\right)$.
Next, let
$$
\mathcal{G}_{8\pi}:=\left\{(\lm,\pl),\,\lm\in(0,8\pi)\right\}\mbox{ be the branch of unique solutions of $\prl$ with $\lm\in(0,8\pi)$.}
$$

In particular, on the basis of the results in \cite{CCL} and \cite{BLin3},
we see that there exists  $\mu>8\pi$ such that
$\mathcal{G}_{8\pi}$ can be continued to a smooth branch,
$$
\mathcal{G}_{\mu}:=\left\{(\lm,\pl),\,\lm\in(0,\mu)\right\},
$$
of solutions of $\prl$ where the first eigenvalue of the linearized problem for $\prl$, which we denote by
$\widehat{\sg}_1(\lm,\pl)$, is strictly positive. Let us define the set of "entropy maximizers",
$$
\mathbb{S}_\om=\bigcup_{E\in [0,+\ii)}\{(\lm,\pl)\,:\;\rl(\pl)\mbox{ is a solution of the \mbox{\bf (MVP)} at energy $E$}\}.
$$

At this point we can explain our motivations. First of all, although it seems obvious in principle,
we don't know of any proof of the fact that the portion of
$\mathcal{G}_{\mu}$ with $\lm>8\pi$ is a subset of $\mathbb{S}_\om$. It looks a rather non trivial task to answer this question in general, since
for a fixed value of $E>E_0$, there can be many solutions of \mbox{\rm $\prl$} with $\lm>8\pi$ such that
\mbox{\rm $\mathcal{E}(\rl(\pl))=E$}, see for example \cite{dem2} for some multiplicity results about Liouville-type equations on compact surfaces.
In the same spirit, the degree evaluation relative to $\prl$, see \cite{CLin2}, \cite{Mal1}, clearly suggests that on non contractible domains the number of solutions
strongly increases with $\lm$. Moreover, on a large class of domains,
one can construct solutions of {\rm $\prl$} attaining large energy values and concentrating at more than one point, see \cite{KMdP}, \cite{EGP},
which however, by {\bf (MVP2)} above, cannot be entropy maximizers for large $E$.
As a matter of fact, it looks hard in general to select the entropy maximizer even among few solutions sharing the same energy.
Actually, this problem is solved in \cite{clmp2} by an argument based on the so called Canonical Variational Principle (CVP for short),
showing that, if $\om$ is of first kind, then $\mathcal{G}_{8\pi}\subset \mathbb{S}_\om$. It is worth to
point out that this method  don't work in case $E>E_{8\pi}$, since the CVP is not well defined for $\lm>8\pi$.
Some results about $\mathcal{G}_{\mu}$ for $\mu$ possibly much larger than $8\pi$ has been recently discussed in \cite{BdM2}.\\
But the situation looks much more involved, since, at least to our knowledge, there are no results at all
in literature about the possible continuation of $\mathcal{G}_{\mu}$, as well as about its relation
with $\mathbb{S}_\om$. For domains of first kind the global connectivity of $\mathcal{G}_{8\pi}$ was proved in \cite{suz}.

\bigskip

Our approach to these problems is based on the following observations. If the energy is monotonic increasing along $\mathcal{G}_{\mu}$, and as
far as $\mathcal{G}_{\mu}$ can be continued to some $\mathcal{G}$ with no bifurcation points, then we can "almost" conclude
that  $\mathcal{G}\subset \mathbb{S}_\om$, see Proposition \ref{enrgm}. After the introduction of a carefully defined set of constrained eigenvalues
($\sg_j=\sg_j(\lm,\pl)$, $j\in\N$) for the linearized
problem relative to $\prl$, see section \ref{sec2} below for further details, we can prove that in fact if $\sg_1>0$, then the energy is monotonic increasing.
In this situation, the implicit function Theorem implies that there are no bifurcation points, whence we can run to the point where
$\sg_1$ vanishes for the "first" time along the branch, say $\lm_*$, where in general we cannot exclude that the branch bifurcates. However, since
the ${\bf (MVP)}$ is a variational maximization problem with two constraints of codimension 1, then we have a natural upper bound for the Morse index.
A rather subtle inspection shows that, under a natural nondegeneracy assumption,
$\mathcal{G}_{\mu}$ can be continued with no bifurcation but only bending points, see Lemma \ref{analytic}, Proposition \ref{prop1} and
Theorem \ref{bifurc}.
This part relies on a careful application of classical techniques \cite{CrRab} together with
the analytic implicit function Theorem \cite{bdt}, \cite{but}. In particular one has to prove that the sign of the variation of
$\lm$ along the branch in a neighborhood of $\lm_*$, is the same as that of $\sg_1$, see Proposition \ref{prcrux}.
The latter information yields the needed monotonicity of the energy and allows
the continuation of $\mathcal{G}_{\lm_*}$ to a branch $\mathcal{G}_{\lm_d}$ of entropy maximizers parametrized by the energy as well.
It is understood that all our results are based on the existence of entropy maximizers, see {\bf (MVP1)}.\\

Let us for the moment just recall that our  modified first eigenvalue,
which we denote by $\sg_1(\lm,\pl)$, satisfies $\sg_1(\lm,\pl)\geq \widehat{\sg}_1(\lm,\pl)$, as far as $\widehat{\sg}_1(\lm,\pl)\geq 0$
(where $\widehat{\sg}_1(\lm,\pl)$ is the standard first eigenvalue). At this point we define,
$$
\lm_*=\lm_*(\om):=\sup\left\{\mu>8\pi\,:\,\sg_1(\lm,\pl)>0,\,\fo\,(\lm,\pl)\in\mathcal{G}_{\mu}\right\},
$$
and, by setting $\mathcal{G}_{\lm_*}:=\left\{(\lm,\pl),\,\lm\in(0,\lm_*)\right\}$,
\beq\label{energy}
E^{(+)}(\lm)=\mathcal{E}\left(\rl(\pl)\right),\quad (\lm,\pl)\in \mathcal{G}_{\lm_*},
\eeq
and
$$
E_*=E_*(\om):=\lim\limits_{\lm\nearrow \lm_*}E^{(+)}(\lm).
$$
Clearly $E^{(+)}(\lm)$ can be extended by continuity on $[0,\lm_*)$ by setting $E^{(+)}(0)=E_0$. Our first result is the following,

\bte\label{thm1}Let $\om$ be a domain of second kind. Then we have,\\
$(i)$ $\lm_*>8\pi$, $E^{(+)}(\lm)$ is analytic, $\frac{d E^{(+)}(\lm)}{d\lm} >0$, $\fo\,\lm\in (0,\lm_*)$ and $E_*>E_{8\pi}$.\\
Moreover, if $\mathbb{S}_\om$ is pathwise connected, it holds:\\
$(ii)$ $\mathcal{G}_{\lm_*}\subset\mathbb{S}_\om$ and there exists a real analytic and strictly increasing function
$\lm^{(+)}_{\mathcal{E}}:(E_0,E_*)\to (0,\lm_*)$ such that
$\mathcal{G}_{\lm_*}=\{(\lm^{(+)}_{\mathcal{E}}(E),\psi_{\sscp \lm^{(+)}_{\mathcal{E}}(E)}),\,E\in (E_0,E_*)\}$ where
$\mbox{\rm $\rh$}_{\sscp \lm^{(+)}_{\mathcal{E}}(E)}(\psi_{\sscp \lm^{(+)}_{\mathcal{E}}(E)})$ is a solution of the
$\mbox{\bf (MVP)}$ at energy $E$ for each $E\in (E_0,E_*)$;\\
$(iii)$ The inverse statistical temperature $\beta(E)=\frac{dS(E)}{d E},\;E\in (E_0,E_*)$, satisfies,
$$
\beta(E)=-\lm^{(+)}_{\mathcal{E}}(E),\;E\in (E_0,E_*),
$$
so that the entropy is strictly decreasing and strictly concave,
$$
\frac{d^2S(E)}{d E^2}=\frac{d\beta(E)}{d E}=-\frac{d\lm^{(+)}_{\mathcal{E}}(E)}{d E}<0,\;E\in (E_0,E_*).
$$
Finally, if we further assume that  $\om$ is strictly starshaped, then,\\
$(iv)$ $E_*<+\ii$, $\lm_*<+\ii$ and $(\lm_*,\psi_{\ssb_*})\in \mathbb{S}_\om$.
\ete

\bigskip

To continue $\mathcal{G}_{\lm_*}$ to a branch of entropy maximizers whose energy $E$ is greater than $E_*$ requires a more subtle argument.
General results yielding the smooth bending of branches of solutions behind the first singular point (i.e. the "first" point where the first eigenvalue vanishes)
for semilinear elliptic equations, see \cite{Amb}, \cite{CrRab}, usually uses the fact that the first eigenvalue is simple and that the first
eigenfunction does not change sign in $\om$. None of these properties is granted in our case, since the linearization of $\prl$ naturally
yields first eigenfunctions which change sign \cite{BLin3}, and first eigenvalues which need not being simple, see for example \cite{B2}.
As mentioned above, we have found a natural physical assumption which allows the solution of these difficulties. In fact, let us first observe that if for some
$E>0$, $\lm=\le$ is the Lagrange multiplier of an entropy maximizer $\rl=\rho_{\le}$, then necessarily we have,
$$
\lm\ino f G[f] -\ino \rl^{-1}f^2\leq 0,\; \fo\,f\in L^{2}(\om)\,:\,\ino f=0,\quad \ino \rl G[f]=\ino\pl f=0,
$$
where the conditions on $f$ are just the compatibility conditions imposed by the mass and energy constraints in the {\bf (MVP)}. Please observe
that this is always true for an entropy maximizer.
The next major improvement of Theorem \ref{thm1} is based on the following definition:

\bdf{\it
Let $E>0$ and $\lm=\le$ be the Lagrange multiplier of an entropy maximizer} $\rl=\rho_{\le}$. {\it We say that}
$\rl$ {\it is a non degenerate entropy maximizer if,}
\beq\label{080216.1}
\lm\ino f G[f] -\ino \rl^{-1}f^2< 0,\; \fo\,f\in L^{2}(\om)\setminus\{0\} \,:\,\ino f=0,\quad \ino \rl G[f]=\ino\pl f=0.
\eeq
\edf

\brm {\it
It is worth to remark that the non degeneracy assumption could be equivalently introduced as the negativity of the first eigenvalue of the quadratic form
in \eqref{080216.1} when restricted to the given vector space. However, the discussion of the continuation
of the branch of entropy maximizers as solutions of} $\prl$ {\it requires a more general spectral analysis, whose
natural ambient space is larger. The key to a major improvement of Theorem \ref{thm1} partially relies on the understanding of the gap
between these two spectral problems, which is the subject of sections \ref{sec2} and \ref{sec4}.}
\erm

By assuming that entropy maximizers are non degenerate in a right neighborhood of $E_*$,
then we obtain the following continuation result for the branch of entropy maximizers $\mathcal{G}_{\sscp \lm_*}$.
\bte\label{thm2}
Let $\om$ be a strictly starshaped domain of second kind and suppose that $\mathbb{S}_\om$ is pathwise connected and that any entropy maximizer with $E\geq E_*$ is non degenerate. Then, there exists $E_d\in (E_*,+\ii)$, $\lm_d\in (8\pi,+\ii)$ and a real analytic function
$\lm_{\mathcal{E}}:(E_0,E_d]\to (0,+\ii)$, such that:\\
$(i)$ $\left.\lm_{\mathcal{E}}\right|_{(E_0,E_*)}\equiv \lm^{(+)}_{\mathcal{E}}$ and $\lm_{\mathcal{E}}(E_d)=\lm_d$;\\
$(ii)$ $\mathcal{G}_{\lm_d}:=\{(\lm_{\mathcal{E}}(E),\psi_{\sscp \lm_{\mathcal{E}}(E)}),\,E\in (E_0,E_d)\}\subset \mathbb{S}_\om$, that is,
for each $E\in (E_0,E_d)$, $\psi_{\sscp \lm_{\mathcal{E}}(E)}$ is a solution of {\rm $\prl$} and 
$\mbox{\rm $\rh$}_{\sscp \lm_{\mathcal{E}}(E)}(\psi_{\sscp \lm_{\mathcal{E}}(E)})$ is a solution of the $\mbox{\bf (MVP)}$
at energy $E$;\\
$(iii)$ there exists $E_m\in [E_*,E_d)$, which is a strict local maximum point of $\lm_{\mathcal{E}}$,
such that $\beta(E)=-\lm_{\mathcal{E}}(E)<0$ in $(E_0,E_d)$, so that $S(E)$ is strictly decreasing in $(E_0,E_d)$,
and either $E_m=E_*$, and then,
\beq\label{convex0}
\graf{\frac{d^2S(E)}{d E^2}=-\frac{d\lm_{\mathcal{E}}(E)}{d E}<  0,\;E\in (E_0,E_*),\\ \\
\frac{d^2S(E)}{d E^2}=-\frac{d\lm_{\mathcal{E}}(E)}{d E}> 0,\;E\in (E_*,E_d),}
\eeq
or $E_m>E_*$ and then $\frac{d^2S(E)}{d E^2}=-\frac{d\lm_{\mathcal{E}}(E)}{d E}\leq 0$ in $(E_0,E_m)$ and
$\frac{d^2S(E)}{d E^2}=-\frac{d\lm_{\mathcal{E}}(E)}{d E}> 0$ in $(E_m,E_d) $, and $\frac{d^2S(E)}{d E^2}=0$ in
$(E_0,E_d)\setminus \{E_m\}$ if and only if $E$ belongs to some finite set $E\in \{E_1,\cdots,E_{n}\}\subset (E_0,E_m)$.
\ete

\bigskip

At least to our knowledge this is the first result about the bending of solutions of $\prl$ behind the first singular point.
The result is also interesting since, by assuming "only" non degeneracy, it yields the existence of an interval (which is $(E_m,E_d)$) of \un{strict convexity} of the Entropy.
Although the statistical temperature in this model has nothing to do with the physical temperature,
it is worth to point out that $(E_m,E_d)$ is an interval of negative (statistical) specific heat, that is, $\frac{d E}{dT_{stat}}<0$.
Actually, if we wish to ensure that $E_*$ is the unique critical point of $\lm_{\mathcal{E}}$, whence in particular the unique absolute maximum point of $\lm_{\mathcal{E}}$,
then, among other things, we need that $\sg_1<0$ for $E>E_*$, (see Theorem \ref{thm3} and fig. 1({\small B}) below)
while Theorem \ref{thm2} would allow in principle for finitely many increasing flex in $(E_*,E_m)$. We think that
this should not be the case for convex domains of second kind, see Conjecture 2 in subsection \ref{hyp}.\\
The proof of Theorem \ref{thm2} uses in a crucial way
the full strength of our constrained eigenvalues and eigenfunctions,
see \rife{lineq0}, \rife{lineq}, which is that they are built up to make the quadratic form in \eqref{080216.1} diagonal in a suitable sense,  see
\eqref{080425.21} and \rife{lineqe} below.\\

It turns out that the situation to go beyond $E_d$ is even more delicate.
By using well known classical results in bifurcation
theory \cite{CrRab}, \cite{Rab} together with those in \cite{bdt}, \cite{but}, then we could follow unambiguously an unbounded branch
of solutions even when higher and higher eigenvalues vanish. Unfortunately neither this fact would suffice to solve our problem.
In fact, even if we would succeed in finding an analytic continuation along $(\lm_*,\psi_{\lm_*})$, it would be
hard in general  to guarantee that the energy is increasing along that path. Actually, this is just
another problem which is hard to control in general, see Remark \ref{rem64} for other details.
Concerning this point, we first state a rather strong result which yields $E_m\equiv E_*$,
$\frac{d^2S(E)}{d E^2}<0$ in $(E_0,E_*)$ and
$\frac{d^2S(E)}{d E^2}>0$ in $(E_*,+\ii)$.
Surprisingly enough, see Remark \ref{stability},  the non degeneracy assumption is still necessary, but it is not sufficient.
It turns out that another natural condition almost does the job.
\bdf{\it
Let $E>0$ and $\lm=\le$ be the Lagrange multiplier of an entropy maximizer} $\rl=\rho_{\le}$. {\it We say that}
$\rl$ {\it is a $\mu_0$-stable entropy maximizer if,}
\beq\label{080216.1nn}
\lm\ino f G[f] -\ino \rl^{-1}f^2\leq \mu_0\ino \rl^{-1}f^2,\; \fo\,f\in L^{2}(\om)\setminus\{0\} \,:\,\ino f=0,\quad \ino \rl G[f]=\ino\pl f=0.
\eeq
\edf

At the end of the day, on strictly starshaped domains of second kind, we
come up with the graph of $S(E)$ as depicted in fig. \ref{fig3} below (which is the same as fig. 4 in \cite{clmp2}).
Here $\sg_i=\sg_i(\lm,\pl)$, $i=1,2$.

\bte\label{thm3}
Let $\om$ be a strictly starshaped domain of second kind and suppose that $\mathbb{S}_\om$ is pathwise connected and that any entropy maximizer with $E\geq E_*$ is non degenerate.
Let $E_d$ be defined as in Theorem \ref{thm2}. If,
for any $(\lm,\pl)\in \{\mathbb{S}_\om \cap \mathcal{E}(\mbox{\rm $\rh$}_{\ssb}(\pl))\geq E_d\}$, it holds $\sg_1<0$ and:\\
${\bf (H1)}$ $\mbox{\rm $\rh$}_{\ssb}$ is $\mu_0$-stable with $\mu_0\leq -\frac{\lm}{\lm+\sg_2}$ and,\\
${\bf (H2)}$ $\sg_2\geq \frac14\frac{\lm |\sg_1|}{\lm +\sg_1}$,\\
then there exists an analytic function $\lm^{(\ii)}_{\mathcal{E}}:(E_0,+\ii)\to (0,+\ii)$, such that:\\
$(i)$ $\left.\lm^{(\ii)}_{\mathcal{E}}\right|_{[E_0,E_d]}\equiv \lm_{\mathcal{E}}$;\\
$(ii)$ $\mathcal{G}_{\ii}:=\{(\lm^{(\ii)}_{\mathcal{E}}(E),\psi_{\sscp \lm^{(\ii)}_{\mathcal{E}}(E)}),\,E\in (E_0,+\ii)\}\equiv \mathbb{S}_\om$, and in
particular, for each $E\in (E_0,+\ii)$, $\psi_{\sscp \lm^{(\ii)}_{\mathcal{E}}(E)}$ is a solution of {\rm $\prl$} and 
$\mbox{\rm $\rh$}_{\sscp \lm^{(\ii)}_{\mathcal{E}}(E)}(\psi_{\sscp \lm^{(\ii)}_{\mathcal{E}}(E)})$ is a solution of the $\mbox{\bf (MVP)}$
at energy $E$;\\\\
$(iii)$ $\lm^{(\ii)}_{\mathcal{E}}(E)\to (8\pi)^+$ as $E\to +\ii$ and  $E_m\equiv E_*$ is the unique critical and maximum point
of $\lm^{(\ii)}_{\mathcal{E}}$. In particular  $\beta(E)=-\lm^{(\ii)}_{\mathcal{E}}(E)$ in $(E_0,+\ii)$ so that $S(E)$ is strictly decreasing in $(E_0,+\ii)$ and
\beq\label{convex}
\graf{\frac{d^2S(E)}{d E^2}=-\frac{d\lm^{(\ii)}_{\mathcal{E}}(E)}{d E}< 0,\;E\in (E_0,E_*),\\ \\
\frac{d^2S(E)}{d E^2}=-\frac{d\lm^{(\ii)}_{\mathcal{E}}(E)}{d E}> 0,\;E\in (E_*,+\ii).}
\eeq
\ete

\bigskip

\brm\label{rem1.11} {\it We remark that in $(ii)$ we have $\mathbb{S}_\om\equiv \mathcal{G}_{\ii}$, that is, under the given assumptions, 
we characterize the full set of entropy maximizers.}
\erm

It will be seen in section \ref{sec2} that $\lm+\sg_1>0$, so that ${\bf (H2)}$ is well posed even if we do not assume the non degeneracy condition.
In particular $\lm+\sg_2>0$ whence ${\bf (H1)}$ is also well posed and stronger than non degeneracy.
See also Remark \ref{stability} for some sufficient conditions which guarantee the $-\frac{\lm}{\lm+\sg_2}$-stability.\\
Whenever the assumptions of Theorem \ref{thm3} are satisfied, then we obtain the global bifurcation diagram of the pairs
$(\lambda^{(\ii)}_{ \mathcal{E}}(E),\psi_{\sscp \lambda^{(\ii)}_{\sscp \mathcal{E}}(E)})$ and
the graph of the entropy $S(E)$, as depicted in figures 1({\small A}), 1({\small B}) and \ref{fig3} respectively.
In particular we have found sufficient conditions which guarantee a positive answer to an open problem about the graph of
$S(E)$, see fig. 4 in \cite{clmp2} and more recently \cite{BdM2} p. 541.  Obviously our result
is consistent with the asymptotic estimates about $S(E)$ for large $E$ as derived in \cite{clmp2}, see Proposition 6.1 in \cite{clmp2} for further details.
We also observe that $S$ is strictly convex in $(E_*,+\ii)$, which is, under our assumptions,
an improvement of Proposition 6.2 in \cite{clmp2}, where it was proved (however with no need of the many hypothesis introduced here)
that $S$ is not concave on domains of second kind.\\
It is interesting to see that, following the route tracked  by mathematical-physics, we come up at once with the solution of a challenging
mathematical problem of independent interest, which is the understanding of the qualitative behavior of the global branch of solutions of $\prl$ emanating
form $\lm=0$.
This was far from trivial already for domains of first kind, where uniqueness and nondegeneracy holds, see \cite{BLin3}, \cite{CCL}, \cite{suz}.
As a matter of fact, and with very few well known exceptions, the qualitative behaviour of global branches of semilinear elliptic problems is poorly understood.
\begin{figure}[h]
    \centering
    \begin{subfigure}[b]{0.3\textwidth}
\psfrag{U}{$<\!\pl\!\!>_{\ssb}$}
\psfrag{V}{${\lambda}$}
\psfrag{L}{${\lambda_*}$}
\psfrag{P}{${8\pi}$}
\psfrag{E}{$2E_*$}
\psfrag{Z}{$2E_0$}
        \includegraphics[width=\textwidth]{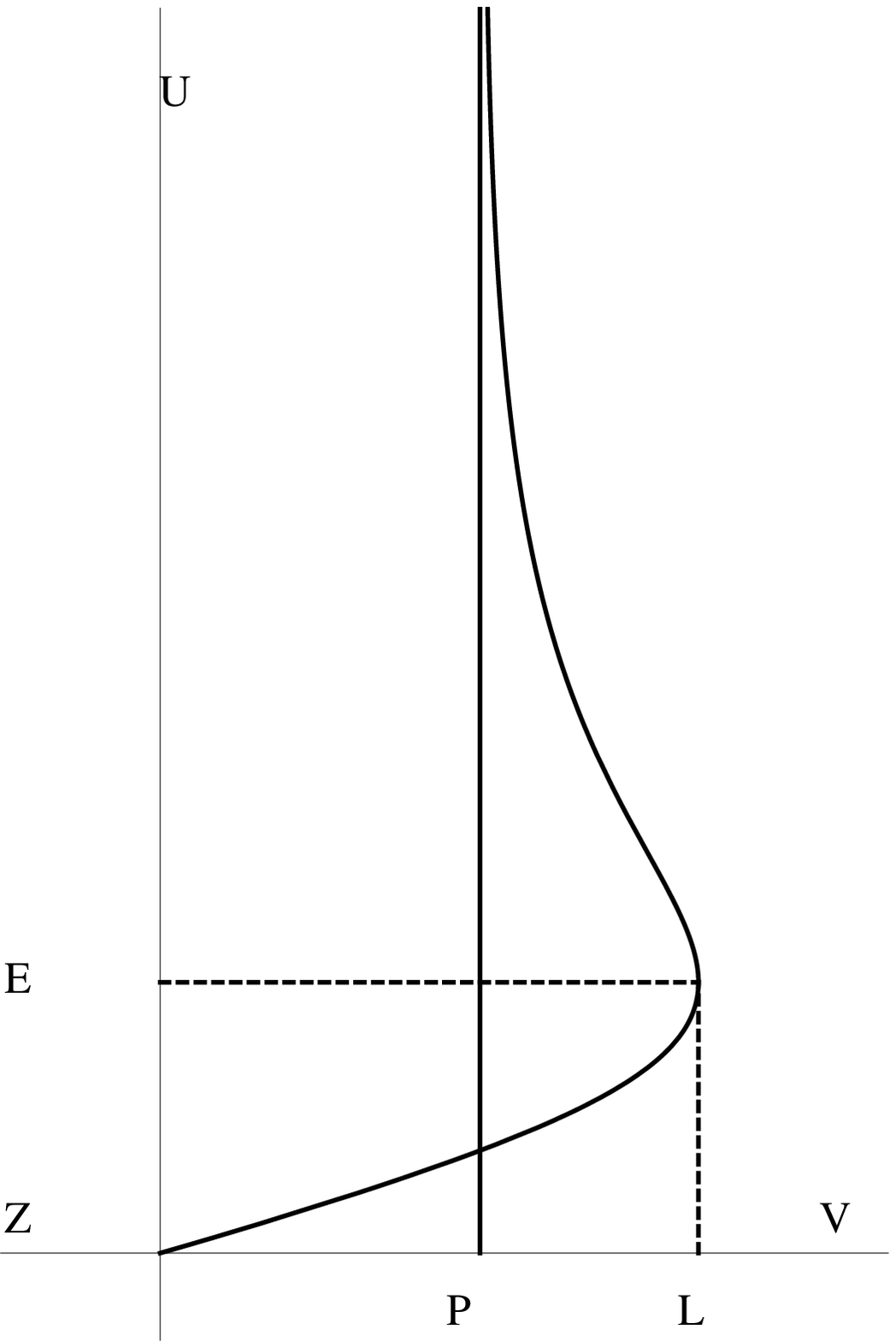}
        \caption{The global branch $\mathcal{G}_{\ii}$ in the plane $(\lm,<\!\pl\!\!>_{\ssb})$.}
        \label{fig1}
    \end{subfigure}
    ~ 
      \hskip4em
    \begin{subfigure}[b]{0.4\textwidth}
    \psfrag{U}{$E$}
\psfrag{V}{${\lambda^{(\ii)}_{ \mathcal{E}}}$}
\psfrag{L}{${\lambda_*}$}
\psfrag{P}{${8\pi}$}
\psfrag{E}{$E_*$}
\psfrag{Z}{$E_0$}
        \includegraphics[width=\textwidth]{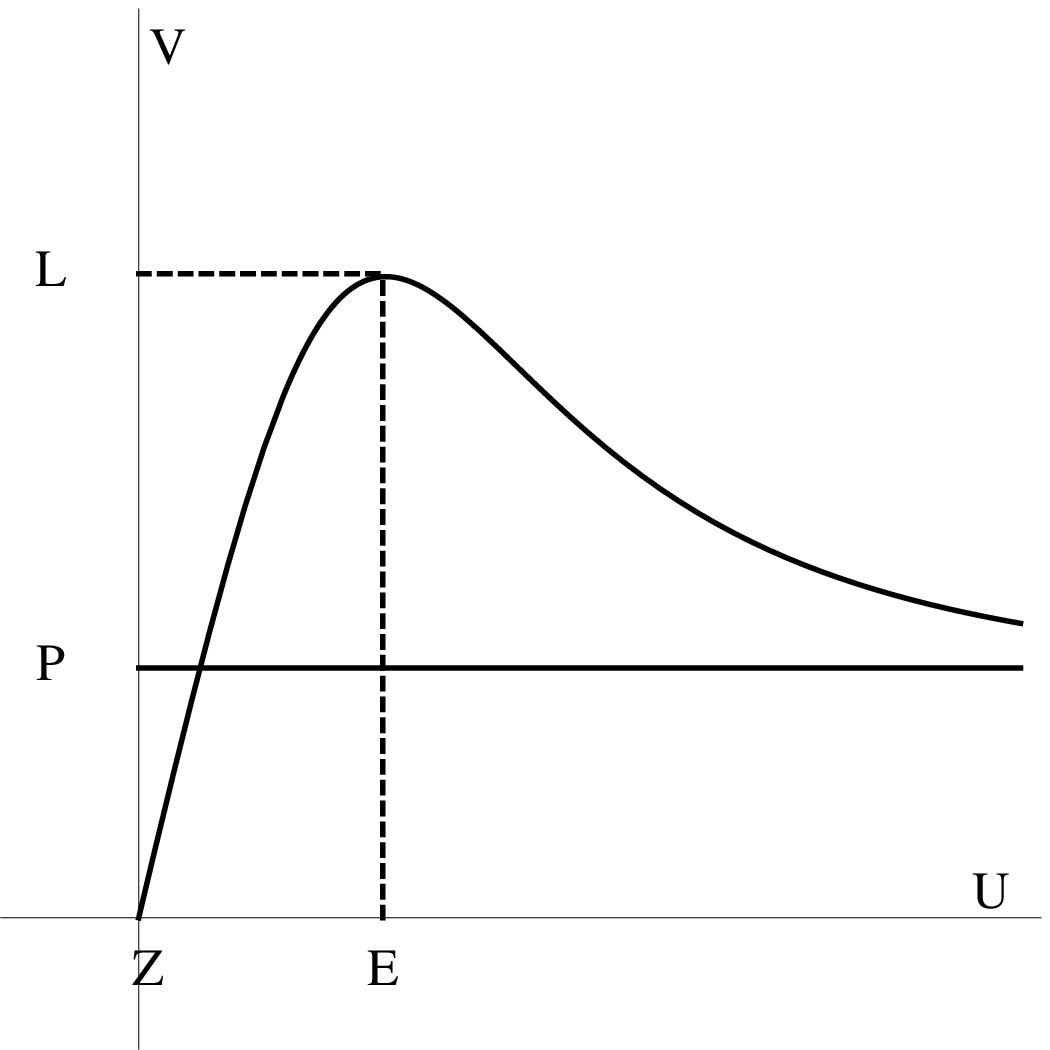}
        \caption{The graph of $\lm^{(\ii)}_{\mathcal{E}}$.}
        \label{fig2}
    \end{subfigure}
    ~ 
     \caption{$\,$}
\end{figure}

\begin{figure}[h]
\psfrag{S}{$S(E)$}
\psfrag{V}{$E$}
\psfrag{M}{$\frac{dS(E)}{d E}\simeq-(8\pi)^+$}
\psfrag{E}{$E_*$}
\psfrag{Z}{$E_0$}
\includegraphics[totalheight=2in]{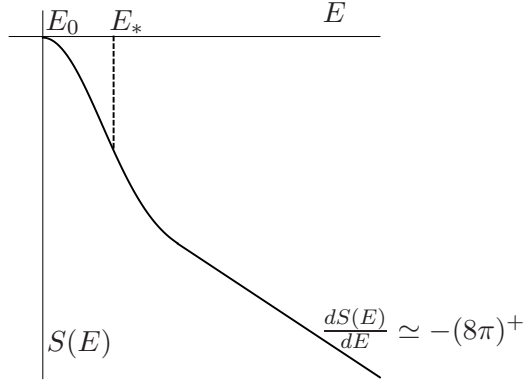}
\caption{The graph of $S(E)$.}\label{fig3}
\end{figure}
In any case, if the assumption about $\sg_1(\lm,\pl)<0$ in Theorem \ref{thm3} is removed, then we still find a global
parametrization of the entropy maximizers, but the entropy has not the nice concave/convex property \rife{convex}, see Theorem \ref{thm4} and
figure \ref{fig4} in section \ref{sec8}.
In fact, in this case
$\lm^{(\ii)}_{\mathcal{E}}$ could have in principle countably many critical points of various sort.
We postpone the discussion of this point to subsection \ref{hyp} below.

\subsection{\bf Conjectures about the global structure of entropy maximizers.}\label{hyp} $\,$\\
It is worth to make some remark about the assumptions made so far.\\
The path-connectedness of $\mathbb{S}_\om$ is very strong from the topological point of view, but otherwise is
rather weak, in the sense that we do not put any other constraint about
the structure of $\mathbb{S}_\om$ and/or its dependence on the values of the energy or of $\lm$. In principle $\mathbb{S}_\om$ could even
be an open set, and one of our goals is to show that it is an analytic curve, see Lemma \ref{enrgm} (which we called "principle of energy parametrization").
Please observe that we are not assuming that the set of solutions of $\prl$ is pathwise connected, which is false in general. In fact,
on any domain of first kind with at least one "hole" \cite{BLin3}, the families of concentrating solutions for $\lm\to 8\pi N$, $N\geq 2$ \cite{KMdP},
\cite{EGP}, are surely not connected with the branch $\mathcal{G}_{8\pi}$. Instead, it seems an interesting open problem
to establish whether or not $\mathbb{S}_\om$ is pathwise connected.\\
However, if we assume that $\mathbb{S}_\om$ is pathwise connected in general, then the non degeneracy does not hold in general.
It has been proved in \cite{CCL} and \cite{BLin3} that if $\gamma=\gamma_\om$ (see {\bf (MVP2)}) has more than one maximum point,
then $\om$ is of second kind. Therefore, in any perfectly symmetric situation of this kind,  the branch of entropy
maximizers would have, for $E>E_*$, many energetically equivalent paths to follow, each one leading to the high energy limit {\bf (MVP3)},
just with distinct concentration points.
Therefore the branch of entropy maximizers soon or later (reasonably at $\lm=\lm_*$) must bifurcate in more than one branch. This is obviously in
contradiction with the non degeneracy property which ensures the lack of bifurcation points.\\
It is reasonable to guess that the entropy maximizers could be proved to be non degenerate, at least
in a suitably defined generic setting relative to some subset of domains whose $\gamma_\om$ have a unique critical point.
This is true for example for convex domains, but it is not true in general for starshaped domains, see \cite{Gu}.
In Theorems \ref{thm2} and \ref{thm3} the domain is assumed to be strictly starshaped,
since in this way we also have that $\lm$ and $E$, are uniformly bounded from above, which is a kind of necessary
condition for all the remaining arguments.\\
It is also reasonable to guess that, for a convex domain of the second kind,  ${\bf (H1)}$ and
${\bf (H2)}$ could be verified at least in the high energy limit, which motivates the following:\\
{\bf Conjecture 1}\\
{\it Let $\om$ be a convex domain of second kind. Then there exists $E_{\om}>E_d$ such that
$\frac{d^2S(E)}{d E^2}=-\frac{d\lm^{(\ii)}_{\mathcal{E}}(E)}{d E}>0$ for any $E>E_{\om}$.}\\

We think at this as a promising route toward the
proof of the strict convexity of the entropy, at least in the high energy limit. We will pursue this approach in a forthcoming paper.\\ 
It is also reasonable to guess that the conclusions of Theorem \ref{thm3} hold for convex domains of second kind,
but this might be far more difficult to prove.\\
{\bf Conjecture 2}\\
{\it Let $\om$ be a convex domain of second kind. Then the conclusions of Theorem \ref{thm3} hold.}\\

In any  case, if $\gamma_\om$ has only one maximum point but also admits
some relative maximum points, then we could face more involved situations. For example, even if the high energy limit would be controlled
by the concentration near the unique maximum point, we could not exclude that the vorticity moves around between those points in the intermediate energy
regime.
So, there could be a chance that entropy maximizers are still non degenerate (so no bifurcation occurs), but
the behavior of $\lm^{(\ii)}_{\mathcal{E}}$ could be more complicated,
with critical points of various sort, which should correspond to the situation depicted in fig. \ref{fig4} rightafter Theorem \ref{thm4}.\\
Although it might look wired in this model, we remark that in higher dimension there are well known examples of branches of semilinear elliptic
equations with exponential and supercritical nonlinearities (arising for example in the study of equilibria of selfgravitating objects \cite{Chandra}, \cite{w})
with infinitely many bending and/or bifurcation points (even in the radial case), see for example \cite{JL} and the more general results in \cite{Dan1}, \cite{Dan2}.

\brm{\it
It is interesting to observe that, with the exception
of the case where $\om$ is a ball, where solutions of \mbox{\rm $\prl$} exist only if $\lm<8\pi$, we still don't know exactly for which
domains of the first kind the set $\mathcal{G}_{8\pi}$ coincides with the set of all solutions of \mbox{\rm $\prl$}.
By assuming Conjecture 2 to be true, the same question could be asked in principle about $\mathcal{G}_{\ii}$ for convex domains of second kind.}
\erm
\brm{\it
Let us remark that, since $\widehat{\sg}_1>0$ along $\mathcal{G}_{\lm_*}$, then one can prove that
if $(\lm,\pl)\in\mathcal{G}_{\lm_*}$, then \mbox{\rm $\rl(\pl)$} is a strict local minimizer of the free energy,
\mbox{\rm $\mathcal{F}_\lm(\rho)=\mathfrak{S}(\rho)+\lm\mathcal{E}(\rho)$}, see \cite{B4}.
Finally, all the result presented in this introduction hold for a more general class of domains, admitting finitely many singular points
of conic type, see \cite{CCL} for rigorous definitions. We skip the details of this facts to avoid technicalities.}
\erm
\brm{\it
The strategy of the proofs of Theorems \ref{thm1}, \ref{thm2}, \ref{thm3} and \ref{thm4} suggests a rather general dual method
to describe the behavior of global branches of solutions for semilinear elliptic problems which seems to deserve further investigations.}
\erm

This paper is organized as follows. In section \ref{sec2} we introduce the "constrained" eigenvalues and eigenfunctions and
prove a version of the analytic implicit function Theorem suitable to be applied to non degenerate entropy maximizers.
In section \ref{sec3} we obtain sufficient conditions which guarantee the monotonicity of the energy. In section \ref{sec4} we discuss
the crucial spectral properties of entropy maximizers. Sections \ref{sec5}, \ref{sec6}, \ref{sec7} are devoted to the proofs of Theorems \ref{thm1},
\ref{thm2}, \ref{thm3}. Section \ref{sec8} is devoted to the proof of Theorem \ref{thm4}.

\section{\bf Preliminary results and notations.}\label{sec2}
Whenever $\pl$ solves $\prl$, and unless no confusion arise, we will denote by $\rl=\rl(\pl)$ the corresponding density.
If $\eta\in  L^2(\om)$  then we set
$$
<\eta>_{\ssb}=\ino \rl \eta \quad\mbox{  and  }\quad [\eta]_{\ol}=\eta \;- <\eta>_{\ssb}.
$$
For the sake of simplicity, with an abuse of notations, we will use the same symbol for functions $\eta_{\ssb}$ already depending on $\lm$, that is,
$$
\mbox{if}\quad <\eta_{\ssb}>_{\ssb}=\ino \rl \eta_{\ssb}\quad\mbox{  then we set }\quad [\eta]_{\ol}=\eta_{\ssb} \;- <\eta_{\ssb}>_{\ssb}.
$$

We also set
$$
<\eta,\varphi>_{\ssb}=\ino \rl \phi\varphi\quad \mbox{and}\quad \|\phi\|_{\ssb}^2=<\phi,\phi>_{\ssb}=\ino \rl \phi^2,
$$
whenever $\{\eta,\varphi\}\subset L^2(\om)$. Since $\rl\in L^{\ii}(\om)$ and it is bounded below away from zero, then
it is easy to see that $<\cdot,\cdot>_{\ssb}$ defines a scalar product on $L^2(\om)$ and that $\|\cdot\|_{\ssb}$ is the corresponding norm.
Next, let us define,
$$
e_0(x)=1,\;x\in \ov{\om},
$$
and for any $\eta\in  L^2(\om)$, possibly depending on $\lm$, denote by,
$$
\al_{\ol}(\eta):=<\eta,e_0>_{\ssb}\equiv <\eta>_{\ssb},
$$
so that $\al_{\ol}(\eta)$ is just the Fourier coefficient of $\eta\in L^2(\om)$ along $X_0:=\mbox{Span}\{e_0\}$, while
$$
[\eta]_{\ol}=\eta-\al_{\ol}(\eta)e_0,
$$
is the projection on $(X_0)^{\perp}$. Moreover, if $\pl$ is a solution of $\prl$, then we will denote by
$$
e_{\ssb,\sse}(x)=\frac{[\psi]_{\ol}}{\|[\psi]_{\ol}\|_{\ssb}},\;x\in \ov{\om},
$$
and similarly, for any $\eta\in L^2(\om)$, define $X_{\ssb,\sse}:=\mbox{Span}\{e_{\ssb,\sse}\}$, and
$$
\al_{\ssb,\sse}(\eta):=<\eta,e_{\ssb,\sse}>_{\ssb},
$$
which is the Fourier coefficient of $\eta$ along $X_{\ssb,\sse}$.
Obviously we have
$$
\|e_{0}\|_{\ssb}=1=\|e_{\ssb,\sse}\|_{\ssb},\;\mbox{and}\;<e_0,e_{\ssb,\sse}>_{\ssb}=0.
$$

\brm\label{rem1}{\it
We will frequently use the fact that, by the Jensen's inequality, we have,
$$
\ino \mbox{\rm $\rl$}[\eta]_{\ol}^2=< [\eta]_{\ol}^2>_{\ssb}=<\eta^2>_{\ssb}-<\eta>_{\ssb}^2\geq 0,
$$
where the equality holds if and only if $\eta$ is constant, whence, in particular, whenever $\eta\in H^{1}_0(\om)$,
if and only if $\eta$ vanishes identically. We will use this property time to time when needed without further comments.}
\erm

\brm {\it
If $\{\eta,\psi\}\subset L^{2}(\om)$, then,
$$
<[\psi]_{\ol},[\eta]_{\ol}>_{\ssb}\equiv <[\psi]_{\ol}, \eta>_{\ssb}\equiv <\psi, [\eta]_{\ol}>_{\ssb},
$$
a fact which will be used time to time when needed, without further comments.}
\erm

For any $\lm\in\R$ let us set,
\beq\label{eF}
F_{\ssb}: C^{2,\alpha}_0(\om)\to  C^{\alpha}(\om),\quad F_{\ssb}(\psi):=-\Delta \psi -\rl(\psi).
\eeq

Next, for any $\pl$ solving $\prl$, we introduce the linearized operator,
\beq\label{2.1}
D_\psi F_{\sscp \lm}(\pl)[\eta]=
-\Delta \eta -\lm \rl [\eta]_{\ol}, \quad \eta \in C^{2,\alpha}_0(\om).
\eeq

We say that $\sg=\sg(\lm,\pl)\in\R$ is an eigenvalue of the linearized operator \rife{2.1} if the equation,
\beq\label{lineq0}
-\Delta \eta-\lm \rl [\eta]_{\ol}=\sg\rl [\eta]_{\ol},
\eeq
admits a non trivial weak solution $\eta\in H^1_0(\om)$. This definition of the eigenvalues requires some comment.
Let $\pl$ be a fixed solution of $\prl$ and let us define,
$$
Y_0:=\left\{ \varphi \in \{L^2(\om),<\cdot,\cdot >_{\ssb}\}\,:\,\ino\rl \varphi\equiv<\varphi,e_0>_{\ssb}=0\right\}.
$$
Clearly $Y_0$ is an Hilbert space, and, since $T(Y_0)\subset W^{2,2}(\om)$, the linear operator,
\beq\label{T0}
T_0:Y_0\to Y_0,\;T_0(\phi)=G[\lm \rl\varphi]-<G[\lm \rl\varphi]>_{\ssb},
\eeq
is self-adjoint and compact. As a consequence,
standard results concerning the spectral decomposition of self-adjoint, compact operators on Hilbert spaces show that
$Y_0$ is the Hilbertian direct sum of the eigenfunctions of $T_0$, which can be represented as $\varphi_k=[\phi_k]_{\ssb,0}$, $k\in\N=\{1,2,\cdots\}$,
$$
Y_0=\overline{\mbox{Span}\left\{[\phi_k]_{\ssb,0},\;k\in \N\right\}}^{L^2(\om)},
$$
for some $\phi_k\in H^1_0(\om)$, $k\in\N=\{1,2,\cdots\}$.
In fact, clearly $\mu=0$ cannot be an eigenvalue of $T_0$ since $T_0(\varphi)=0$ readily implies $\varphi\equiv 0$.
Therefore, the eigenfunction $\varphi_k$, whose eigenvalue is $\mu_k=\frac{\lm}{\lm+\sg_k}\in \R\setminus\{0\}$, satisfies,
$$
\varphi_k=(\lm+\sg_k) \left(G[\rl\varphi_k]-< G[\rl\varphi_k]>_{\ssb}\right),
$$
that is, by defining,
$$
\phi_k:=(\lm+\sg_k) G[\rl\varphi_k],
$$
it is easy to see that $\varphi_k$ is an eigenfunction of $T_0$ if and only if
$\phi_k$ is in $H^1_0(\om)$ and weakly solves,
\beq\label{lineq}
-\Delta \phi_k= \tau_k \rl [\phi_k]_{\ol}\quad\mbox{ in }\quad \om,\quad \tau_k=\lm+\sg_k.
\eeq

In particular we will use over and over the fact that $\varphi_k=[\phi_k]_{\ssb,0}$ and
\beq\label{lineq9}
\phi_k=(\lm+\sg_k) G[\rl[\phi_k]_{\ssb,0}],\quad k\in\N=\{1,2,\cdots\}.
\eeq

At this point, standard arguments in the calculus of variations show that,
\beq\label{4.1}
\sg_1=\sg_1(\lm,\pl)=\inf\limits_{\phi \in H^1_0(\om)\setminus \{0\}}
\dfrac{\ino |\nabla \phi|^2 - \lm \ino \rl[\phi]_{\ssb,0}^2 }{\ino  \rl[\phi]_{\ssb,0}^2}.
\eeq
The ratio in the right hand side of  \rife{4.1} is well defined because of
Remark \ref{rem1}. Higher eigenvalues are defined inductively via the variational problems,
\beq\label{4.1.k}
\sg_k=\sg_k(\lm,\pl)=\inf\limits_{\phi \in H^1_0(\om)\setminus \{0\},<\phi,[\phi_m]_{\ssb,0}>_{\ssb}=0,\; m\in \{1,\ldots, k-1\}}
\dfrac{\ino |\nabla \phi|^2 - \lm \ino \rl[\phi]_{\ssb,0}^2 }{\ino  \rl[\phi]_{\ssb,0}^2}.
\eeq
Obviously, if $\phi_i$ and $\phi_j$ are eigenfunctions corresponding to distinct eigenvalues $\sg_i\neq \sg_j$, then
\beq\label{orto}
<[\phi_i]_{\ssb,0}, [\phi_j]_{\ssb,0}>_{\ssb}=0,
\eeq
as also easily follows from \rife{lineq}. The eigenvalues form a numerable non decreasing
sequence $\sg_1(\lm,\pl)\leq \sg_2(\lm,\pl)\leq....\leq \sg_k(\lm,\pl)\leq...\,$, where in particular,
\beq\label{lineq2}
\tau_k=\lm+\sg_k\geq \lm +\sg_1>0,\;\fo\,k\in \N,
\eeq
the last inequality being an immediate consequence of \rife{4.1}. Obviously, by the Fredholm alternative, see \cite{GL},
if $0\notin \{\sg_j\}_{j\in \N}$, then $I-T_0$ is an isomorphism of $Y_0$ onto itself.\\
Finally, any $\psi \in L^{2}(\om)$ admits the Fourier series expansion,
\beq\label{Four}
\psi=\al_{\lm,0}(\psi)+\sum\limits_{j=1}^{+\ii} \al_{\lm,j}(\psi)[\phi_{j}]_{\lm,0},\quad \|[\phi_{j}]_{\lm,0}\|_{\ssb}=1,
\eeq
where
$$
\al_{\lm,j}(\psi)=<[\phi_{j}]_{\lm,0},\psi>_{\ssb}.
$$

If we let $\widehat{\sg}_1$ be the standard first eigenvalue defined by,
$$
\widehat{\sg}_1=\inf\limits_{\phi \in H^1_0(\om)\setminus \{0\}}\dfrac{\ino |\nabla \phi|^2 - \lm \ino \rl(\phi^2-<\phi>_{\ssb}^2) }{\ino  \rl\phi^2},
$$
then we see that, as far as $\widehat{\sg}_1(\lm,\pl)\geq 0$, we have,
\beq\label{larger}
\sg_1\geq \widehat{\sg}_1,
\eeq
and then, in view of the results in \cite{suz}, as later improved in \cite{CCL} and \cite{BLin3}, we have the following:\\

{\bf Theorem A}[\cite{BLin3}, \cite{CCL}, \cite{suz}]
{\it Let $\om$ be a domain of second kind. For any $(\lm,\pl)\in \ov{\mathcal{G}_{8\pi}}$, the first eigenvalue of
\eqref{2.1} is strictly positive, that is,  $\sg_1(\lm,\pl)>0$, for each $\lm\in[0,8\pi]$.}

\bigskip

The following observation will be crucial in the sequel. It states that the projections of any eigenfunction along $e_0$ and
$e_{\ssb,\mathcal{E}}$ are proportional.
\bpr
Let $\pl$ be any solution of {\rm $\prl$, $\rl=\rl(\pl)$} be the corresponding density and $\phi_k$ be any eigenfunction of \eqref{2.1} with
eigenvalue $\sg_k$. Then
\beq\label{pro.16.1}
\al_{\ol}(\phi_k)=\|[\psi]_{\ol}\|_{\ssb}\tau_k \al_{\ssb,\sse}(\phi_k).
\eeq
\epr
\proof
Integrating \rife{lineq} we see that $\ino\Delta \phi_k=0$, whence we easily conclude that
$$
-\ino(\Delta \phi_k)[\psi]_{\ol}=-\ino(\Delta \phi_k)\pl.
$$

Therefore we find
$$
\|[\psi]_{\ol}\|_{\ssb}\tau_k \al_{\ssb,\sse}(\phi_k)=\tau_k\ino \rl \phi_k [\psi]_{\ol}=\tau_k\ino \rl [\phi_{k}]_{\ol} [\psi]_{\ol}=-\ino(\Delta \phi_k)[\psi]_{\ol}=
$$
$$
-\ino(\Delta \phi_k)\pl=-\ino\phi_k\Delta\pl=\ino \rl \phi_k= \al_{\ssb,0}(\phi_k),
$$
where the penultimate equality follows  integrating by parts.
\finedim

\bigskip

Finally, we will frequently use the following Lemma.
The first part is a standard application of the analytic implicit function Theorem, see \cite{Dan} and
more recently \cite{but}, to our problem. The second part is a careful adaptation of well known arguments, see \cite{Amb}, which yields
a detailed description of the invertibility of $D_{\psi}F_{\ssb}(\psi_{\ssb})$ when $0$
is a simple eigenvalue.
\ble\label{analytic} Let $\psi_{\sscp \lm_0}$ be a solution of {\rm $\prl$} with $\lm=\lm_0$.\\
If $0$ is not an eigenvalue of \eqref{2.1}. Then:\\
$(i)$ $D_{\psi}F_{\sscp \lm_0}(\psi_{\sscp \lm_0})$ is an isomorphism of $C^{2,\al}_0(\om)$ onto $C^\al(\om)$;\\
$(ii)$ There exists a neighborhood $\mathcal{U}\subset
\R\times C^{2,\al}_0(\om)$ of $(\lm_0,\psi_{\sscp \lm_0})$ such that the set of solutions of
{\rm $\prl$} in $\mathcal{U}$ is an analytic curve $I\ni\lm\mapsto \pl\in B$, for
suitable neighborhoods $I$ of $\lm_0$ in $\R$ and $B$ of $\psi_{\sscp \lm_0}$ in $C^{2,\al}_0(\om)$.\\
If $0$ is a simple eigenvalue of \eqref{2.1}, that is, its eigenspace  $\mbox{\rm Span}\{\phi_1\}$ is one dimensional, then\\
$(iii)$ $D_{\psi}F_{\sscp \lm_0}(\psi_{\sscp \lm_0})$ is an isomorphism of $Y_1$ onto $R$, where,
$$
Y_1=\left\{\xi \in C^{2,\al}_{0}(\om)\,:\,\ino \mbox{\rm $\rh_{\sscp \lm_0}$} [\phi_{1}]_{\sscp \lm_0,0}\xi =0\right\}\mbox{ and }
R=\left\{h\in C^{\al}(\om)\,:\, \ino\phi_1 h=0.\right\}.
$$

\ele
\proof
It is not difficult to check that $F_{\ssb}$ is jointly analytic \cite{but}, with respect to $(\lm,\psi)\in \R\times C^{2,\al}_0(\om)$.
Therefore, whenever $(i)$ holds, then $(ii)$ is an immediate consequence of the analytic implicit function theorem, see \cite{Dan} and more recently \cite{but}.\\
In the rest of the proof we set $\lm=\lm_0$.\\
$(i)$  Let $h\in C^{\al}(\om)$, then we have to prove that the equation,
$$
-\Delta \phi -\lm \rl [\phi]_{\ssb,0}=h \;\; \mbox{ in }\;\;\om,
$$
admits a unique solution in $C^{2,\al}_0(\om)$, or either, by standard elliptic regularity theory \cite{GL}, equivalently,
that the equation,
\beq\label{un1}
\phi=G[\lm\rl[\phi]_{\ssb,0}]+G[h],
\eeq
admits a unique solution $\phi_h \in C^\al(\om)$. Let $X_{0,\al}=X_0\cap C^{\al}(\om)$ and $Y_{0,\al}=Y_0\cap C^{\al}(\om)$.
Since $X_{0,\al} \oplus Y_{0,\al} = C^\al(\om)$, after projection on $X_{0,\al}$ and $Y_{0,\al}$ we see that
\rife{un1} is equivalent to the system,
$$
\graf{[\phi]_{\ssb,0}= [G[\lm\rl[\phi]_{\ssb,0}]+G[h]]_{\ssb,0},\\ \\<\phi>_{\ssb}=<G[\lm\rl[\phi]_{\ssb,0}]+G[h]>_{\ssb}.}
$$
Since $0\notin \{\sg_j\}_{j\in\N}$, then $I-T_{0}$ is an isomorphism of $Y_{0,\al}$ onto itself. Therefore, the first equation, which has the form $(I-T_0)([\phi]_{\ssb,0})=[G[h]]_{\ssb,0}$,
has a unique solution $\varphi_h\in Y_{0,\al}$. Let
$$
\phi_{h}:=G[\lm\rl\varphi_h]+G[h]=\varphi_h +<G[\lm\rl\varphi_h]+G[h]>_{\ssb}.
$$
Clearly $\phi_h\in C_0^{2,\al}(\om)$ and in particular it is the unique function which satisfies both
the equations.\\

$(iii)$ If $h\in R$, then $h_\lm=(\rl)^{-1}h\in C^{\al}(\om)$ satisfies $\ino \rl \phi_1h_\lm =0$ and we can
prove equivalently, as in $(i)$, that the equation,
\beq\label{un2}
\phi=G[\lm\rl[\phi]_{\ssb,0}]+G[\rl h_\lm],\mbox{ with the condition }\ino \rl \phi_1h_\lm =0,
\eeq
admits a unique solution $\phi_h\in C^{\al}(\om)$ such that $\ino \rl[\phi_{1}]_{\ol} \phi_h =0$.\\
Let $\phi_{1,0}=[\phi_1]_{\ol}$, $X_{1,0,\al}=\mbox{Span}\{\phi_{1,0}\}\cap C^{\al}(\om)$ and
$Y_{1,0,\al}=(X_{0,\al}\oplus X_{1,0,\al})^{\perp}\cap C^{\al}(\om)$.
Note that, putting $\phi=\al_0+\al_1\phi_{1,0}+\varphi$ and $h_{\lm}=\al_{0,h}+\al_{1,h}\phi_{1,0}+f_h$, where
$\al_0=\al_{\ssb,0}(\phi)$, $\al_1=\al_{\ssb,1}(\phi)$, $\al_{0,h}=\al_{\ssb,0}(h)$,  $\al_{1,h}=\al_{\ssb,1}(h)$,
$\{\varphi,f_h\}\in Y_{1,0,\al}$, and in view of \rife{lineq9} with $\sg_1=0$, then \rife{un2} takes the  form,
$$
\al_0+\varphi=G[\lm\rl \varphi]+G[\rl f_{h}]+\al_{0,h}\pl+\frac{\al_{1,h}}{\lm}\phi_{1,0}+\left(\al_1+\frac{\al_{1,h}}{\lm}\right)\al_0(\phi_{1}),
$$
with the conditions
$$
\al_{1,h}=-\al_{0,h}\al_{0}(\phi_1) \iff \ino \rl \phi_1h_\lm =0,
$$
and
$$
\al_1=0 \iff \ino \rl [\phi_{1}]_{\ol}\phi=0.
$$
By using the Fourier decomposition \rife{Four} together with \rife{lineq9}, we conclude that,\\
$\{[G[\lm\rl \varphi]]_{\ol}, [G[\rl f_h]]_{\ol} \}\in Y_{1,0,\al}$.
At this point, setting $\al_{1,\psi}=\al_{\ssb,1}(\pl)$, we can project along $X_{0,\al}$, $X_{1,0,\al}$ and $Y_{1,0,\al}$, to obtain the system,

$$
\graf{\varphi=[G[\lm\rl \varphi]]_{\ol}+[G[\rl f_h]]_{\ol}+\al_{0,h}([\pl]_{\ol}-\al_{1,\psi}\phi_{1,0}),\\ \\
0= \al_{0,h}\al_{1,\psi}+\frac{\al_{1,h}}{\lm},\\ \\
\al_0=<G[\lm\rl \varphi]>_{\ssb}+<G[\rl f_{h}]>_{\ssb}+\,\al_{0,h}<\pl>_{\ssb}+\left(\al_1+\frac{\al_{1,h}}{\lm}\right)\al_0(\phi_{1}),\\ \\
\al_{1,h}=-\al_{0,h}\al_{0}(\phi_1),\\ \\
\al_1=0.}
$$

\bigskip

Plugging the fourth into the second we obtain,

$$
\al_{0,h}\left(\al_{1,\psi}-\frac{\al_{0}(\phi_1)}{\lm}\right)=0.
$$
This equation is always satisfied regardless of the value of $\al_{0,h}$ as it is a particular case of \rife{pro.16.1},
$$
\al_{0}(\phi_1)=\ino \rl \phi_1=\ino (-\Delta \pl) \phi_1=\ino \pl (-\Delta \phi_1)=\lm\ino \rl \phi_{1,0}\pl=\lm \al_{1,\psi}.
$$
Therefore, \rife{un2} is equivalent to the reduced system,

$$
\graf{\varphi=[G[\lm\rl \varphi]]_{\ol}+[G[\rl f_h]]_{\ol}+\al_{0,h}([\pl]_{\ol}-\al_{1,\psi}\phi_{1,0}),\\ \\
\al_0=<G[\lm\rl \varphi]>_{\ssb}+<G[\rl f_{h}]>_{\ssb}+\,\al_{0,h}<\pl>_{\ssb}+\frac{\al_{1,h}}{\lm}\al_0(\phi_{1}).
}
$$
The first equation has the form
\beq\label{un3}
\left.(I-T_0)\right|_{Y_{1,0,\al}}(\varphi)=\eta,\quad \eta \in Y_{1,0,\al},
\eeq

and since the homogeneous equation $\varphi=\left.T_0\right|_{Y_{1,0,\al}}\varphi$ has no non trivial solutions, then by the Fredholm alternative, \rife{un3}
admits a unique solution $\varphi_h\in Y_{1,0,\al}$. Clearly $\varphi_h$ uniquely defines $\al_{0}$,
whence the solution $\phi_h=\al_0+\varphi_h$ is unique and satisfies $\phi_h\in C^{\al}(\om)$, as claimed.
\finedim

\bigskip
\bigskip

\section{\bf Monotonicity of the energy along $\mathcal{G}_{\lm_*}$.}\label{sec3}

Our first result is about the monotonicity of the energy for solutions of $\prl$.

\bpr\label{pr2.1} Let $\pl$ be a solution of {\rm $\prl$}.\\
$(i)$ Whenever the function,
$$
\vb=\dfrac{d \pl}{d \lm},
$$
is pointwise well defined in $\om$,  then we have,
{\rm
\beq\label{8.12.11}
\frac{d}{d \lm} \mathcal{E}\left(\rl(\psi_\lm)\right)=\left<\vb \right>_{\ssb}=<[\psi]_{\ol}^2>_{\ssb}+\lm<[\psi]_{\ol} v_{\ssb}>_{\ssb}.
\eeq
}

$(ii)$ {\rm If $\sg_1(\lm,\pl)>0$ then $<[\psi]_{\ol} v_{\ssb}>_{\ssb}> 0$.
In particular, if $\vb$ is pointwise well defined and if either $\sg_1(\lm,\pl)\geq 0$  or  $<\psi_{\ol} v_{\ol}>_{\ssb}\geq 0$, then $\frac{d}{d \lm} \mathcal{E}\left(\rl(\psi_\lm)\right)>0$.}

\bigskip

$(iii)$ Let $\om$ be a domain of second kind and $E^{(+)}(\lm)$ be defined in \eqref{energy}. Then $E^{(+)}(\lm)$ is analytic and satisfies,
{\rm
\beq\label{8.12.1}
\frac{d}{d \lm} E^{(+)}(\lm)=\frac{d}{d \lm} \mathcal{E}\left(\rl(\psi_\lm)\right)>0,\;\fo\;\lm\in(0,\lm_*).
\eeq}
\noi \!\!\!\! In particular $E^{(+)}(\lm)$ is strictly increasing in $(0,\lm_*)$
and defines an analytic bijection which maps  $(0,\lm_*)$ onto $(E_0,E_*)$. Clearly its inverse $\lm^{(+)}(E)$
is analytic and strictly increasing and satisfies $\lm^{(+)}(E)\searrow 0$, as $E\searrow E_0$, $\lm^{(+)}(E_{8\pi})=8\pi$
and $\lm^{(+)}(E)\nearrow \lm_*$ as $E\nearrow E_*$.
\epr

\proof $(i)$ If $\pl$ is differentiable as a function of $\lm$, then standard elliptic arguments based on $\prl$ and the associated Green's representation formula, show that
$\vb\in C^{2,\alpha}_0(\om)$ is a classical solution of,

\beq\label{1b1}
-\Delta \vb =\rl[\psi]_{\ol}+\lm \rl [v]_{\ol}.
\eeq
In particular,  since the right hand side of $\prl$ is analytic as a function of $\lm$, then, by continuous dependence on the data,
we see that $\pl$ is analytic as a function of $\lm$ in some small enough open interval. At this point, by using $\prl$ we also find,

$$
\mathcal{E}\left(\rl\right)=\frac12\ino \rl\pl=\frac12\ino |\nabla \pl|^2,
$$

and then we conclude that,

\beq\label{8.12.2}
\frac{d}{d \lm} \mathcal{E}\left(\rl\right)=\ino (\nabla \vb,\nabla \pl)=-\ino \vb (\Delta \pl)=<\vb>_{\ssb}.
\eeq
Next, by using $\prl$ and \eqref{1b1}, we see that
$$
<\vb>_{\ssb}=\ino\rl\vb =\ino -(\Delta \pl)\vb=\ino -\pl(\Delta \vb)=
$$
\beq\label{2b1}
<[\psi]_{\ol}^2>_{\ssb}+\lm<[\psi]_{\ol} [v]_{\ol}>_{\ssb},
\eeq
which in view of \rife{8.12.2} proves \rife{8.12.11}.\\

$(ii)$ If $\sg_1(\lm,\pl)>0$, then $\vb$ is well defined by Lemma \ref{analytic} and we can use \rife{8.12.11}.
Let,
$$
[\pl]_{\lm,0}=\sum\limits_{j=1}^{+\ii}\al_j\phi_{j,0}, \quad
[\vb]_{\lm,0}=\sum\limits_{j=1}^{+\ii}\beta_j\phi_{j,0},
$$
be the Fourier expansions \rife{Four} of $[\pl]_{\lm,0}$ and $[\vb]_{\lm,0}$ with
respect to the normalized, $\|\phi_{j,0}\|_{\lm}=1$, projections $\phi_{j,0}:=[\phi_j]_{\lm,0}$.
By using \rife{1b1} with \rife{lineq9}, we obtain,
\beq\label{lamq31}
\sg_{j}\ino \rl\phi_{j,0}[\vb]_{\lm,0} =
\ino\rl\phi_{j,0}[\pl]_{\lm,0},\mbox{ that is }\beta_j=\frac{\al_j}{\sg_{j}},
\eeq
where $\sg_{j}=\sg_{j}(\lm,\pl)$. As a consequence, we find,
\beq\label{event1}
<[\pl]_{\lm,0},[\vb]_{\lm,0}>_{\lm}=\sum\limits_{j=1}^{+\ii} \al_j\beta_j=
\sum\limits_{j=1}^{+\ii} \sg_{j}(\beta_j)^2\geq \sg_{1}<[\vb]_{\lm,0}^2>_{\lm}> 0,
\eeq
whenever $\sg_1(\lm,\pl)> 0$, which proves the first part of $(ii)$. Next, by virtue of \rife{8.12.11} and Remark \ref{rem1},
we see that  $\frac{d}{d \lm} \mathcal{E}\left(\rl(\psi_\lm)\right)><[\pl]_{\lm,0},[\vb]_{\lm,0}>_{\lm}$. As a consequence,
by \rife{event1} we conclude that
$\frac{d}{d \lm} \mathcal{E}\left(\rl(\psi_\lm)\right)>0$ whenever $\vb$ is well defined and
either $\sg_1(\lm,\pl)\geq 0$  or  $<[\pl]_{\ol} [\vb]_{\ol}>_{\ssb}\geq 0$, as claimed.

\bigskip
\bigskip

$(iii)$ Since $\sg_1(\lm,\pl)>0$ for any $(\lm,\pl)\in\mathcal{G}_{\lm_*}$, then Lemma \ref{analytic} implies
that $\mathcal{G}_{\lm_*}$ is an analytic curve parametrized by $\lm$ and in particular that $\pl$ is an analytic function of $\lm$. Since $\mathcal{E}$
is also jointly analytic with respect to $\lm$ and $\psi$, then
$E^{(+)}(\lm)$ is analytic as well, and we can apply $(i)-(ii)$ to conclude that \rife{8.12.1} holds.
The rest of the statement is an obvious consequence of the monotonicity of $E^{(+)}(\lm)$.
\finedim

\bigskip
\bigskip

\section{\bf Spectral properties of Entropy maximizers.}\label{sec4}
The analysis of the spectral properties of entropy maximizers is based on the following equivalent formulation of \rife{080216.1}.
For fixed $\lm=\le>0$ and if $\rl=\rl(\pl)$ is an entropy maximizer, then $\pl=\ple$ solves $\prl$ and  $\rl\in C^{0}(\ov{\om})$ and $\rl>0$ on $\ov{\om}$.
Therefore, a function $f\in L^{2}(\om)$ satisfies $\ino f=0$ and $ \ino \rl G[f]\equiv\ino\pl f=0$ if and only if it takes the form,
\beq\label{080216.11}
f=\rl \varphi,\quad \varphi\in \mathbb{T}_{\ssb}:=\left\{ \varphi\in L^{2}(\om)\,:\, \quad \ino \rl \varphi=0,\quad \ino\rl  \pl \varphi=0\right\}.
\eeq
Clearly a function $\varphi_{\sscp T}\in L^2(\om)$ satisfies $\varphi_{\sscp T}\in \mathbb{T}_{\ssb}$ if and only if,
$$
\varphi_{\sscp T}=\varphi-\al_{\ol}(\varphi)e_0-\al_{\ssb,\sse}(\varphi)e_{\ssb,\sse}=[\varphi]_{\ol}-\al_{\ssb,\sse}(\varphi)e_{\ssb,\sse},
$$
for some $\varphi\in L^2(\om)$. In fact, we obviously have $\ino \rl \varphi_{\sscp T}=0$, and,
$$
\ino \rl \pl \varphi_{\sscp T}=\ino \rl\pl [\varphi]_{\ol}-<[\varphi]_{\ol},e_{\ssb,\sse}>_{\ssb}\ino \rl \pl e_{\ssb,\sse}=
$$
$$
\ino  \rl\pl[\varphi]_{\ol}-<[\varphi]_{\ol},\frac{[\psi]_{\ol}}{\|[\psi]_{\ol}\|_\lm^2}>_{\ssb}\ino  \rl\pl [\psi]_{\ol}\equiv
\ino \rl\pl [\varphi]_{\ol}-
<\pl,[\varphi]_{\ol}>_{\ssb}=0.
$$

Therefore we can substitute $f=\rl\varphi_{\sscp T}$ in \rife{080216.1} to obtain,

\beq\label{080425.21}
\mathcal{A}_{\lm}(\varphi_{\sscp T})<0,\; \varphi_{\sscp T}=[\varphi]_{\ol}-\al [\psi]_{\ol},\;
\fo\,\varphi\in L^{2}(\om)\,:\;\varphi \notin  X_0\oplus X_{\ssb,\sse},
\eeq

where $\mathcal{A}_{\lm}$ is the quadratic form
\beq\label{080425.21nn}
\mathcal{A}_{\lm}(\varphi)=\lm\ino \rl\varphi G[\rl\varphi] -\ino \rl\varphi^2,\quad
\fo\,\varphi\in L^{2}(\om),
\eeq

and $\al$ is defined as follows,
$$
\al=\dfrac{\al_{\ssb,\sse}(\varphi)}{\|[\psi]_{\ol}\|_\lm}.
$$

Similarly, when we will need the $\mu_0$-stability property \rife{080216.1nn}, then we will use the following inequality,

\beq\label{080216.1nnn}
\mathcal{A}_{\lm}(\varphi_{\sscp T})\leq \mu_0 \ino \rl \varphi^2_{\sscp T},\; \varphi_{\sscp T}=[\varphi]_{\ol}-\al [\psi]_{\ol},\;
\fo\,\varphi\in L^{2}(\om)\,:\;\varphi \notin  X_0\oplus X_{\ssb,\sse}.
\eeq
We will use several times the following fact: if, for some $\lm>0$, $\phi_k$ is an eigenfunction of \rife{2.1} with eigenvalue $\sg_k$, then
the following holds,

\beq\label{lineqe}
\mathcal{A}_{\lm}([\phi_{k}]_{\ol})=\lm\ino \rl[\phi_{k}]_{\ol} G[\rl[\phi_{k}]_{\ol}] -\ino \rl[\phi_{k}]_{\ol}^2=-\frac{\sg_k}{\tau_k}\ino \rl[\phi_{k}]_{\ol}^2,
\eeq
where we used \rife{lineq9} and $\tau_k=\lm+\sg_k>0$, see \rife{lineq2}.

\bigskip

\bpr\label{prop1}
For fixed $E>0$ let $\lm=\le$ be a Lagrange multiplier of a non degenerate entropy maximizer {\rm $\rl$} and
$\pl$ be the corresponding solution of {\rm $\prl$}. Then:\\
$(i)$ If $\phi_k$ is any eigenfunction of \eqref{2.1} corresponding to an
eigenvalue $\sg_k\leq 0$, then $\al_{\ssb,0}(\phi_k)\neq 0\neq \al_{\ssb,\mathcal{E}}(\phi_k)$;\\
$(ii)$   \eqref{2.1} admits at most one non positive eigenvalue;\\
$(iii)$ If $\sg_k\leq 0$ is a non positive eigenvalue of \eqref{2.1}, then it is simple, that is, it admits only one eigenfunction.
\epr
\proof
$(i)$ Arguing by contradiction, if we would have $\al_{\ol}(\phi_k)=0$, then, in view of \rife{pro.16.1}, we would also have $\al_{\ssb,\sse}(\phi_k)=0$ and
so we would conclude that $\al_{\ssb,\sse}(\phi_k)= 0$. Consequently, we could test \rife{080425.21} with $\varphi_{\sscp T}\equiv \phi_k$, and then
we would find,
$$
\mathcal{A}_\lm(\phi_k)=\lm\ino \rl\phi_k G[\rl\phi_k]-\ino \rl\phi^2_k<0.
$$
However this is impossible since then, by using \rife{lineqe}, we would also obtain,
$$
0>\lm\ino \rl\phi_k G[\rl\phi_k]-\ino \rl\phi^2_k=\frac{|\sg_k|}{\tau_k}\ino \rl\phi^2_k,
$$
which is the desired contradiction. Then, in view of \rife{lineq2} and \rife{pro.16.1}, we see that $\al_{\ssb,\mathcal{E}}(\phi_k)\neq0$ as well.
\bigskip

$(ii)$ We argue by contradiction and suppose that there exist two eigenfunctions $\phi_i$, $i=1,2$ corresponding to a pair of distinct eigenvalues satisfying
$\sg_1< \sg_2\leq 0$. In view of  $(i)$ we have $ \al_{\ol}(\phi_i)\neq 0\neq\al_{\ssb,\sse}(\phi_i)$, $i=1,2$ and so we define,
$$
\phi_{t_1,t_2}=t_1(\phi_1-\al_{\ol}(\phi_1)e_0 -\al_{\ssb,\sse}(\phi_1)e_{\ssb,\sse})+
t_2(\phi_2-\al_{\ol}(\phi_2)e_0-\al_{\ssb,\sse}(\phi_2)e_{\ssb,\sse})\in \mathbb{T}_{\ssb}.
$$
Clearly, choosing
$$
t_1=t_{1,0}:=\al_{\ssb,\sse}(\phi_2),\quad t_2=t_{2,0}:=-\al_{\ssb,\sse}(\phi_1),
$$
we also have
$$
\phi_{t_{1,0},t_{2,0}}=t_{1,0}\phi_{1,0}+t_{2,0}\phi_{2,0},\;\mbox{where}\;\phi_{i,0}=[\phi_{i}]_{\ol},\,i=1,2.
$$

Therefore we can test \rife{080425.21} with $\varphi_{\sscp T}=\phi_{t_{1,0},t_{2,0}}$, to conclude that
$$
\mathcal{A}_{\lm}(\phi_{t_{1,0},t_{2,0}})=\lm\ino \rl \phi_{t_{1,0},t_{2,0}}G[\rl \phi_{t_{1,0},t_{2,0}}] -\ino \rl\phi_{t_{1,0},t_{2,0}}^2< 0.
$$

Since $\sg_1\neq \sg_2$, then \rife{orto} implies that,
$$
\ino\rl\phi_{1,0}\phi_{2,0}=\ino\rl\phi_{1,0}\phi_{2}=0=\ino\rl\phi_{1}\phi_{2,0}=\ino\rl\phi_{1,0}\phi_{2,0},
$$
and so, by using \rife{lineqe}, we also find,
$$
0>\mathcal{A}_{\lm}(\phi_{t_{1,0},t_{2,0}})=
$$
$$
t_{1,0}^2\lm \ino\rl\phi_{1,0}G[\rl\phi_{1,0}] + t_{2,0}^2\lm \ino\rl\phi_{2,0}G[\rl\phi_{2,0}] +2\lm t_{1,0}t_{2,0} \ino\rl\phi_{1,0}G[\rl\phi_{2,0}]
$$
$$
-t_{1,0}^2\ino\rl\phi_{1,0}^2-t_{2,0}^2\ino\rl\phi_{2,0}^2-2t_{1,0}t_{2,0}\ino\rl\phi_{1,0}\phi_{2,0}=$$
$$
\left(\lm \frac{t_{1,0}^2}{\tau_1}-t_{1,0}^2\right)\ino\rl\phi_{1,0}^2+\left(\lm \frac{t_{2,0}^2}{\tau_2}-t_{2,0}^2\right)\ino\rl\phi_{2,0}^2+2\lm\frac{t_{1,0}t_{2,0}}{\tau_2}\ino\rl\phi_{1,0}\phi_{2,0}=
$$
$$
\frac{t_{1,0}^2}{\tau_1}|\sg_1|\ino\rl\phi_{1,0}^2+\frac{t_{2,0}^2}{\tau_2}|\sg_2|\ino\rl\phi_{2,0}^2,
$$
which is the desired contradiction.\\

$(iii)$ We argue by contradiction and suppose that there exist two eigenfunctions $\phi_i$, $i=1,2$ corresponding to the eigenvalue
$\sg_1\leq 0$. In view of  $(i)$ we have $ \al_{\ol}(\phi_i)\neq 0\neq\al_{\ssb,\sse}(\phi_i)$, $i=1,2$, and so we define,
$$
\phi_{t_1,t_2}=t_1(\phi_1-\al_{\ol}(\phi_1)e_0 -\al_{\ssb,\sse}(\phi_1)e_{\ssb,\sse})+
t_2(\phi_2-\al_{\ol}(\phi_2)e_0-\al_{\ssb,\sse}(\phi_2)e_{\ssb,\sse})\in \mathbb{T}_{\lm}.
$$
This time, by choosing,
$$
t_1=t_{1,0}:=\al_{\ssb,\sse}(\phi_2),\quad t_2=t_{2,0}:=-\al_{\ssb,\sse}(\phi_1),
$$
we find,
$$
\phi_{t_{1,0},t_{2,0}}=t_{1,0}\phi_{1}+t_{2,0}\phi_{2}.
$$

Therefore we can can test \rife{080425.21} with $\varphi_{\sscp T}=\phi_{t_{1,0},t_{2,0}}$, to conclude that
$$
\mathcal{A}_{\lm}(\phi_{t_{1,0},t_{2,0}})=\lm\ino \rl \phi_{t_{1,0},t_{2,0}}G[\rl \phi_{t_{1,0},t_{2,0}}] -\ino \rl\phi_{t_{1,0},t_{2,0}}^2< 0.
$$

In the same time, by using \rife{lineqe}, we find,
$$
0>\mathcal{A}_{\lm}(\phi_{t_{1,0},t_{2,0}})=
$$
$$
t_{1,0}^2\lm \ino\rl\phi_{1}G[\rl\phi_{1}] + t_{2,0}^2\lm \ino\rl\phi_{2}G[\rl\phi_{2}] +2\lm t_{1,0}t_{2,0} \ino\rl\phi_{1}G[\rl\phi_{2}]
$$
$$
-t_{1,0}^2\ino\rl\phi_{1}^2-t_{2,0}^2\ino\rl\phi_{2}^2-2t_{1,0}t_{2,0}\ino\rl\phi_{1}\phi_{2}=$$
$$
\left(\lm \frac{t_{1,0}^2}{\tau_1}-t_{1,0}^2\right)\ino\rl\phi_{1}^2+\left(\lm \frac{t_{2,0}^2}{\tau_1}-t_{2,0}^2\right)\ino\rl\phi_{2}^2+
2\left( \lm\frac{t_{1,0}t_{2,0}}{\tau_1}-t_{1,0}t_{2,0}\right)\ino\rl\phi_{1}\phi_{2}=
$$
$$
\frac{t_{1,0}^2}{\tau_1}|\sg_1|\ino\rl\phi_{1}^2+\frac{t_{2,0}^2}{\tau_1}|\sg_1|\ino\rl\phi_{2}^2+
2\frac{t_{1,0}t_{2,0}}{\tau_1}|\sg_1|\ino\rl\phi_{1}\phi_{2}=
$$
$$
\frac{|\sg_1|}{\tau_1}\ino \rl\phi_{t_{1,0},t_{2,0}}^2
$$
which is the desired contradiction.\finedim

\bigskip
\bigskip

\section{\bf The proof of Theorem \ref{thm1}.}\label{sec5}

In this section we  prove Theorem \ref{thm1}.

\bigskip

{\bf The proof of Theorem \ref{thm1}.}\\

$(i)$ In view of Theorem A in section \ref{sec2},  if $\om$ is of second kind, then $\sg_1(8\pi,\psi_{\sscp 8\pi})>0$. Therefore, by Lemma \ref{analytic}, we can
continue $\mathcal{G}_{8\pi}$ in a small  right neighbourhood of $\lm_{8\pi}$ where $\sg_1(\lm,\pl)>0$. As a consequence $\lm_*>8\pi$, as claimed.
The rest of the proof of $(i)$ is an immediate consequence of Proposition \ref{pr2.1}.\\

$(ii)$ We first prove that $\mathcal{G}_{\lm_*}\subset \mathbb{S}_\om$. Actually we prove a slightly more general result which will be needed later on.

\bpr[Principle of Energy parametrization]\label{enrgm}$\,$\\
Let $\om$ be a domain of second kind and suppose that $\mathbb{S}_{\om}$ is pathwise connected.\\
If ${\Gamma}=\{(\lm(t),\psi_{\ssb(t)}),\; t\in (0,t_0)\}$ satisfies:\\
$(a)$ for each $t\in (0,t_0)$, $\psi_{\ssb(t)}$ is a solution of {\rm $\prl$} with $\lm=\lm(t)$ and $\ov{\mathcal{G}_{8\pi}}\subset \Gamma$;\\
$(b)$ the map $t\to (\lm(t),\psi_{\ssb(t)})$ is injective and analytic, $\lm(t)\to 0^+$ as $t\to 0^+$, $\lm(t_c)=8\pi$ for some $t_c\in (0,t_0)$ and
$\lm^{'}(t)>0$ for $t\in (0,t_c]$;\\
$(c)$ $\Gamma$ has no bifurcation points, that is, for any $(\lm_1,\psi_{\lm_1})\in\Gamma$ the set of solutions of {\rm $\prl$}
in a small enough neighbourhood $\mathcal{U}\in \R\times C^{2,\al}_0(\om)$ of $(\lm_1,\psi_{\lm_1})$ coincides with
$\mathcal{U}\cap \Gamma$;\\
$(d)$ the energy satisfies $\frac{d}{dt}\mathcal{E}(\mbox{\rm $\rh$}_{\ssb(t)}(\psi_{\ssb(t)}))>0$ in $(0,t_0)$.\\
Then, the inclusion $\Gamma\subset \mathbb{S}_{\om}$ holds and there exists an analytic function
$\lm^{(0)}_{\mathcal{E}}:(E_0,E_{t_0})\to (0,+\ii)$, where $E_{t_0}=\lim\limits_{t\to t_0} \mathcal{E}(\mbox{\rm $\rh$}_{\ssb(t)}(\psi_{\ssb(t)}))$,
such that
$\Gamma=\{(\lm^{(0)}_{\mathcal{E}}(E),\psi_{\sscp \lm^{(0)}_{\mathcal{E}}(E)}),\,E\in (E_0,E_{t_0})\}$ where
$\mbox{\rm $\rh$}_{\sscp \lm^{(0)}_{\mathcal{E}}(E)}(\psi_{\sscp \lm^{(0)}_{\mathcal{E}}(E)})$ is a solution of the $\mbox{\bf (MVP)}$ at energy $E$
for any $E\in (E_0,E_{t_0})$.
\epr

\proof
We argue by contradiction and suppose that there exists $(\lm_s,\psi_{\sscp \lm_s})\in \Gamma$ such that $\lm_s=\lm(t_s)$ for some $t_s\in (0,t_0)$
and
$(\lm_s,\psi_{\sscp \lm_s})\notin \mathbb{S}_\om$. Since $\ov{\mathcal{G}_{8\pi}}\subset\mathbb{S}_\om$, then by $(b)$ we find
$\lm_s>8\pi$, $t_s>t_c$ and we can
set $\Gamma_s=\{(\lm(t),\psi_{\ssb(t)}),\; t\in (0,t_s)\}$ and $E_s=\mathcal{E}(\rho_{\sscp \lm_s}(\psi_{\sscp \lm_s}))$.
In view of $(d)$, the energy is strictly increasing along $\Gamma$, whence
$E_s>\mathcal{E}(\mbox{\rm $\rh$}_{\ssb(t_c)}(\psi_{\ssb(t_c)}))= E_{8\pi}$, with $E_{8\pi}$ as in Definition \ref{def2}.
In particular, if we fix $E_2=\mathcal{E}(\rho_{\sscp \lm_2}(\psi_{\sscp \lm_2}))>E_s$, where $\rho_{\sscp \lm_2}(\psi_{\sscp \lm_2})\in \mathbb{S}_\om$
(which always exists by {\bf (MVP1)}), and since
$\mathbb{S}_\om$ is pathwise connected by assumption and $\ov{\mathcal{G}_{8\pi}}\subset\mathbb{S}_\om$, then there must
exist a continuous path $\{(\lm_x,\psi_{\sscp \lm_x}),x\in [1,2]\} \subset \mathbb{S}_\om$ such that
$(\lm_1,\psi_{\sscp \lm_1})\in\ov{\mathcal{G}_{8\pi}}$, for some $\lm_1\leq 8\pi$ (so that in particular
$(\lm_1,\psi_{\sscp \lm_1})\in\Gamma$, since $\ov{\mathcal{G}_{8\pi}}\subset \Gamma$ by $(a)$) while
$(\lm_2,\psi_{\sscp \lm_2})\in \mathbb{S}_\om\setminus \Gamma_s$, since by $(d)$ we have
$\sup\limits_{(\lm,\pl)\in\Gamma_s}\mathcal{E}(\rl(\pl))=E_s$. In other words,
we see that for $x$ close enough to $2$, we have $(\lm_x,\psi_{\sscp \lm_x})\notin \Gamma_s$,
while $(\lm_x,\psi_{\sscp \lm_x})\in \Gamma_s$ for some $x\geq 1$. Therefore, since
$\{(\lm_x,\psi_{\sscp \lm_x}),x\in [1,2]\}$ is continuous, then it is well defined,
$$
\widehat{x}=\inf\{x\in[1,2]\;:\;(\lm_x,\psi_{\sscp \lm_x})\notin \Gamma_s\},
$$
which satisfies $\widehat{x}\in[1,2)$ and $(\lm_{\widehat{x}},\psi_{\sscp \lm_{\widehat{x}}})\notin\Gamma_s$.
Indeed, since each $(\lm_x,\psi_{\sscp \lm_x})$ is a solution of $\prl$, then
$(\lm_{\widehat{x}},\psi_{\sscp \lm_{\widehat{x}}})$ would be a bifurcation point on $\Gamma$, in the sense that the set of solutions of $\prl$
in any small enough neighborhood $\mathcal{U}\in \R\times C^{2,\al}_0(\om)$ of $(\lm_{\widehat{x}},\psi_{\sscp \lm_{\widehat{x}}})$ would
not be just $\mathcal{U}\cap \Gamma$. This is in contradiction with $(c)$, and then we find
$(\lm_{\widehat{x}},\psi_{\sscp \lm_{\widehat{x}}})\notin\Gamma_s$. Since $(\lm_{\widehat{x}},\psi_{\sscp \lm_{\widehat{x}}})\neq
(\lm_s,\psi_{\sscp \lm_s})$ by assumption, then necessarily $\lm_{\widehat{x}}=0$ and then $(\lm_x,\psi_{\sscp \lm_x})$ has non empty intersection
with the set of solutions of $\prl$ for some $\lm\leq 0$. This is neither possible since the set of solutions of $\prl$ for $\lm\in(-\ii,8\pi)$ is a smooth curve of entropy maximizers whose energy is less than $E_0$ and
but with no bifurcations points. Therefore $\Gamma\subset \mathbb{S}_\om$, as claimed.\\
Next observe that $E^{(0)}(t):=\mathcal{E}(\mbox{\rm $\rh$}_{\ssb(t)}(\psi_{\ssb(t)}))$ is analytic in $(0,t_0)$, and by $(d)$ it satisfies,
$(E^{(0)})^{'}(t)>0$ in $(0,t_0)$. Therefore it admits an analytic inverse, see \cite{but}, which we denote by $\mathcal{T}:(E_0,E_{t_0})\to (0,t_0)$,
with the property that for each $t\in(0,t_0)$, there exists a unique $E\in(E_0,E_{t_0})$ such that
$\lm^{(0)}_{\mathcal{E}}(E):=\lm(\mathcal{T}(E))\equiv \lm(t)$. Obviously  $\lm^{(0)}_{\mathcal{E}}$ is analytic, as claimed.
\finedim

\bigskip
\bigskip

{\bf The proof of Theorem \ref{thm1} continued.}\\
Since $\sg_1(\lm,\pl)>0$ whenever $(\lm,\pl)\in\mathcal{G}_{\lm_*}$, then, by Lemma  \ref{analytic},
we see that $\mathcal{G}_{\lm_*}$ is an analytic curve with no bifurcation points.
In view of Proposition \ref{pr2.1}, $\mathcal{E}(\rl(\pl))$ is monotonic increasing along $\mathcal{G}_{\lm_*}$.
Then Proposition \ref{enrgm} can be applied and we let $\lm^{(+)}_{\mathcal{E}}$ to be the (analytic) inverse of $E^{(+)}$ as
defined in Proposition \ref{pr2.1}. Therefore,
$\mathcal{G}_{\lm_*}=\{(\lm^{(+)}_{\mathcal{E}}(E),\psi_{\sscp \lm^{(+)}_{\mathcal{E}}(E)}),\,E\in (E_0,E_*)\}$ where
$\mbox{\rm $\rh$}_{\sscp \lm^{(+)}_{\mathcal{E}}(E)}(\psi_{\sscp \lm^{(+)}_{\mathcal{E}}(E)})$ is a solution of the $\mbox{\bf (MVP)}$ at energy $E$,
and the proof of $(ii)$ is complete.\\

\bigskip

$(iii)$ To simplify the notations let us write,
$$
\lm=\lm^{(+)}_{\mathcal{E}}(E),\;
\lm^{'}:=\frac{d\lm^{(+)}_{\mathcal{E}}(E)}{dE},\; \pl=\psi_{\lm^{(+)}_{\mathcal{E}}(E)},\;
\vb=\left.\frac{d\pl}{d\lm}\right|_{\lm=\lm^{(+)}_{\mathcal{E}}(E)}.
$$
By using $\prl$ we find,
$$
S(E)=-2\lm E+\log\left(\ino e^{\lm\pl}\right),
$$
and then obtain,
\beq\label{formula0}
\frac{dS(E)}{dE}=-2E\lm^{'}-2\lm +\lm^{'}{\ino \rl \pl}+\lm\lm^{'}{\ino \rl \vb}=-2\lm+\lm\frac{1}{\ino \rl \vb}\ino \rl \vb=-\lm,
\eeq
where we used \rife{8.12.11}, which proves the first part of the claim. In particular we find
\beq\label{formula}
\frac{d^2S(E)}{dE^2}=-\lm^{'},
\eeq
which is always negative by $(ii)$ and so $(iii)$ is proved as well.

\bigskip

\brm\label{formulas} {\it Obviously the proof of $(iii)$ shows that \eqref{formula0} and \eqref{formula} hold,
whenever $S(E)$ is defined and twice differentiable on an energy value $E\in (E_0,+\ii)$.}
\erm
\bigskip

$(iv)$ If $\om$ is strictly starshaped, then by the Pohozaev identity, see for example \cite{clmp1}, we know that there exists $\widehat{\lm}\geq 8\pi$,
depending only by $\om$, such that
\beq\label{poho}
\sup\left\{\lm>0\,:\,\pl\mbox{ is a solution of }\prl \right\}\leq \widehat{\lm}.
\eeq

Therefore $\lm_*<+\ii$ and in particular $\lm_*>8\pi$, since $\om$ is of second kind. We conclude the proof with a Lemma which will
be needed later on.
\ble\label{ebound} If $\om$ is of second kind and $(\lm(E),\psi_{\lm(E)})\in \mathbb{S}_\om$ for a
continuous function $\lm(E)$, and if $\lm_n=\lm(E_n)\to \lm_*>8\pi$,
as $E_n\to E_*$, then $E_*<+\ii$.\\
In particular $\psi_{\sscp \lm_n}\to \psi_*$ uniformly, where $\psi_*$ is a solution of {\rm $\prl$} with $\lm=\lm_*$, and $(\lm_*,\psi_{\lm_*})\in \mathbb{S}_\om$.
\ele
\proof
We assume by contradiction that $E_*=+\ii$. Since
$$
2E=2\mathcal{E}(\mbox{\rm $\rh_{\sscp \lm(E)}$}(\psi_{\sscp \lm(E)}))=<\psi_{\sscp \lm(E)}>_{\sscp \lm(E)},
$$
then there exists a subsequence (which we won't relabel) $\lm_n\to \lm_*$, such that $u_n:=\lm_n\psi_{\lm_n}$ satisfies  $\|u_n\|_{\ii}\to +\ii$ as $n\to +\ii$.
Since $u_n\in C^{2,\al}_0(\om)$ satisfies the Liouville type equation $-\Delta u_n=\lm_n\dfrac{e^{u_n}}{\ino e^{u_n}}$ in $\om$, then
by Lemma 2.1 in \cite{CCL} we know that $u_n$ is uniformly bounded in a small enough neighbourhood of
$\pa\om$. Therefore $u_n$ admits finitely many blow up points in $\om$ in the sense of Brezis-Merle \cite{bm}.
More exactly, by well known concentration-compactness results \cite{yy}, we conclude that,
$$
\lm_n\rh_n=\lm_n\dfrac{e^{u_n}}{\ino e^{u_n}}\rightharpoonup 8\pi \sum\limits_{j=1}^{m}\delta_{p_j},\ainf,
$$
weakly in the sense of measures in $\om$, where $\{p_1,\cdots,p_m\}\subset \om$ and $\delta_p$ is the Dirac delta with singular point $p$. However $\rh_n$ is a sequence of entropy maximizers whose energy is
$E_n=\mathcal{E}(\mbox{\rm $\rh_n$}(u_n))\to +\ii$, and then, according to {\bf (MVP2)}, it has one concentration point, that is $m=1$. Therefore we find,
$$
8\pi =\lim\limits_{n\to+\ii} \lm_n\ino \rh_n=\lm_n\to \lm_*>8\pi,
$$
which is a contradiction. Therefore $E_*<+\ii$ and, as a consequence, by the same argument we find $\|\pl\|_\ii\leq C$, as $\lm\to \lm_*$,
for a suitable constant $C>0$, see also Lemma 2.1 in \cite{BdM2}.
At this point, by a standard bootstrap argument, we have that $\psi_{\sscp \lm_n}\to\psi_*$, as $E_n\to E_*$, where
$\psi_{*}$ is a solution of $\prl$ with $\lm=\lm_*$. Obviously, since $S(E)$ is continuous (see {\bf (MVP2)}) and since
$\rh_{\sscp \lm_*}(\psi_*)$ is the uniform limit of entropy maximizers, then it is
an entropy maximizer as well, that is $(\lm_*,\psi_{\lm_*})\in \mathbb{S}_\om$, as claimed. \finedim
\bigskip

Clearly Lemma \ref{ebound} yields the desired conclusion and so the proof of $(iv)$ is completed.
\finedim

\bigskip
\bigskip

\section{\bf The proof of Theorem \ref{thm2}.}\label{sec6}
In this section we prove Theorem \ref{thm2}. To achieve our goal we will need the following continuation result for non degenerate entropy maximizers.
We will use some facts about analytic maps in Banach spaces and the analytic implicit function theorem, and we refer to \cite{but}
for the details about these well known results.
\bte\label{bifurc}
Let $\widehat{\lm}>0$ and $\widehat{\psi}:=\psi_{\sscp \widehat{\lm}}$ satisfy $(\widehat{\lm},\widehat{\psi})\in\mathbb{S}_\om$
and suppose that $\sg_1(\widehat{\lm},\widehat{\psi})=0$. If $\mbox{\rm $\rl$}(\pl)$ is non degenerate, then $\sg_1(\widehat{\lm},\widehat{\psi})$
is simple, that is, the kernel
of \eqref{2.1} at $\lm=\widehat{\lm}$ is one dimensional, say $X_{1}=\mbox{\rm Span}\{\phi_1\}$,
where $\|[\phi_1]_{\widehat{\lm},0}\|_{\widehat{\lm}}=1$. Let
$$
Y_1=\left\{u\in C^{2,\alpha}_0(\om)\,:\, \ino u \phi_1=0\right\}.
$$
Then there exists $\eps>0$ and $\lm:(-\eps,\eps) \to \R$, $u:(-\eps,+\eps)\to Y_1$,
which are analytic functions of $s$, and satisfy $\lm(0)=\widehat{\lm}$, $u(0)=0$ and,
\beq\label{210216.1.n}
\psi(s):=\widehat{\psi}+s\phi_1+ u(s),\quad \mbox{which satisfies}\quad \psi(0)=\widehat{\psi},
\eeq
is a solution of
\beq\label{3.1.n}
\graf{-\Delta \psi(s) =\mbox{\rm \rh}_{\sscp \lm(s)}(\psi(s)) \quad \om, \\ \psi(s) = 0 \qquad \qquad\pa\om,}
\eeq

for any $s\in(-\eps,\eps)$. Putting $\frac{d}{ds}:={}^{'}$, the first order derivatives of $\lm(s),u(s)$ satisfy,
$$
u^{'}(0)=0,\quad \psi^{'}(0)=\phi_1,\quad \lm^{'}(0)=0.
$$
\ete
\proof
Since $(\widehat{\lm},\widehat{\psi})\in\mathbb{S}_\om$ then
Proposition \ref{prop1} $(ii)-(iii)$ ensure that the kernel of the
linearized operator \rife{2.1} at $\lm=\widehat{\lm}$ is
one dimensional, which is denoted by $X_1=\mbox{span}\{\phi_1\}$. Let
$\rh_z=\rh_{\widehat{\lm}}(\widehat{\psi}\,)$. Then we have,
\beq\label{3.1.1.n}
D_\psi F_{\sscp \widehat{\lm}}(\widehat{\psi}\,)[\phi_1]=-\Delta \phi_1-\widehat{\lm} \rh_z[\phi_1]_{\widehat{\lm},0}=0.
\eeq

Therefore, the map $\Phi:\R\times \R \times Y_1\to C^{\alpha}(\om)$, defined by,
$$
\Phi (s,\mu,u)=-\Delta (\widehat{\psi} + s \phi_1 + u)- \rh_{\sscp \widehat{\lm} +\mu}(\widehat{\psi} + s \phi_1 + u),
$$
is jointly analytic with respect to its variables, see \cite{but}, it satisfies,
$$
\Phi (0,0,0)=-\Delta\widehat{\psi} -\rh_z=0,
$$
and its derivative $D_{\mu,u}\Phi (0,0,0):\R \times Y_1\to C^{\alpha}(\om)$,
$$
D_{\mu,u}\Phi (0,0,0)=\left(-\rh_z[\widehat{\psi}\,]_{\sscp  \widehat{\lm},0} ,\;D_\psi F_{\sscp  \widehat{\lm}_*}(\widehat{\psi}\,) \right),
$$

acts on $\R \times Y_1$ as follows,
$$
\R \times Y_1\ni (\mu, u)\mapsto D_{\mu,u}\Phi (0,0,0) [\mu,u]=
-\rh_z[\widehat{\psi}\,]_{\sscp  \widehat{\lm},0}\mu+D_\psi F_{\sscp  \widehat{\lm}}( \widehat{\psi}\,)[u].
$$

By Lemma \ref{analytic} $(iii)$ we find that
$D_\psi F_{\sscp  \widehat{\lm}}( \widehat{\psi}\,)$ is an isomorphism of $Y_1$ onto its range $R$, that is,
$$
R=\left\{h\in C^{\al}(\om)\,:\,\ino h \phi_1=0 \right\}.
$$

We claim that,
\ble
$D_{\mu,u}\Phi (0,0,0)$ is an isomorphism if and only if
\beq\label{crux.n}
\widehat{a}:= {\|[\widehat{\psi}\,]_{\sscp \widehat{\lm},0}\|_{\widehat{\lm}}}\al_{\widehat{\lm},\sse}(\phi_1)\equiv
 \ino \mbox{\rm $\rh$}_z[\widehat{\psi}\,]_{\sscp \widehat{\lm},0}\phi_1=<[\widehat{\psi}\,]_{\sscp \widehat{\lm},0},\phi_1>_{\sscp \widehat{\lm}}\neq 0.
\eeq
\ele
\proof
Since $X_{1}\oplus R=C^{\al}(\om)$, then the equation $D_{\mu,u}\Phi (0,0,0) [\mu,u]=h\in C^{\al}(\om)$ is equivalent to the system
$$
\graf{-\mu \ino \rh_z[\widehat{\psi}\,]_{\sscp \widehat{\lm},0}\phi_1= \ino h \phi_1, \\ \\
P_R(-\rh_z[\widehat{\psi}\,]_{\sscp \widehat{\lm},0}) \mu + D_\psi F_{\sscp \widehat{\lm}}(\widehat{\psi}\,) u = P_R(h), }
$$

where $P_R$ is the projection onto $R$. Obviously, the first equation is solvable if and only if $\widehat{a}\neq 0$, while the second equation is equivalent to
$$
D_\psi F_{\sscp \widehat{\lm}}(\widehat{\psi}) u=P_R(h+\mu\rh_z[\widehat{\psi}\,]_{\sscp \widehat{\lm},0}),
$$
which always admits a unique solution, since $D_\psi F_{\sscp \widehat{\lm}}(\widehat{\psi}\,)$ is an isomorphism of $Y_1$ onto its range $R$.\finedim

\bigskip

At this point we just observe that, in view of Proposition \ref{prop1} $(i)$, \rife{crux.n} is surely satisfied, whence
the analytic implicit function theorem can be applied, see for example \cite{but}, to conclude that in a small interval
$s\in(-\eps,\eps)$ there exist
$\lm:(-\eps,\eps) \to \R$ and $u:(-\eps,+\eps)\to Y_1$, which are analytic functions of $s$, satisfying,
$\lm(0)=\widehat{\lm}$, $u(0)=0$, and in particular that \eqref{210216.1.n} and \eqref{3.1.n} holds.
The first order derivatives of $(\lm(s),u(s))$ as a function of $s$, are linear in
$$
D_s\Phi (0,0,0)\equiv D_\psi F_{\sscp \lm^*}(\psi_*)[\phi^*_1]=0,
$$
whence, putting $\frac{d}{ds}:={}^{'}$, we find,
$$
u^{'}(0)=0,\quad \psi^{'}(0)=\phi^*_1,\quad \lm^{'}(0)=0,
$$
as claimed.\finedim

\bigskip
\bigskip

At this point we are ready to present,\\
{\bf The proof of Theorem \ref{thm2}.}\\
Since $\om$ is strictly starshaped and of second kind, and since $\mathbb{S}_\om$ is assumed to be pathwise connected, then Theorem \ref{thm1}
applies and we have $\lm_*<+\ii$, $E_*<+\ii$ and $(\lm_*,\psi_*)\in \mathbb{S}_\om$, where $\psi_*:=\psi_{\sscp \lm_*}$. Clearly
$\sg_1(\lm_*,\psi_*)=0$ and, since $\rh_{\sscp \lm_*}$ is non degenerate by assumption, then we can apply Theorem \ref{bifurc}, to conclude that
the kernel of the
linearized operator \rife{2.1} at $\lm=\lm_*$ is
one dimensional, and we will denote it by $X_{\lm_*,1}=\mbox{span}\{\phi^*_1\}$, where $\|[\phi^*_1]_{\lm_*,0}\|_{\lm_*}=1$.
Let $\rh_*=\rh_{\sscp \lm_*}\equiv\rh_{\sscp \lm_*}(\psi_*)$ and let us define,
$$
Y_*=\left\{u\in C^{2,\alpha}_0(\om)\,:\, \ino u \phi^*_1=0\right\}.
$$
Then, by Theorem \ref{bifurc}, there exist
$\lm:(-\eps,\eps) \to \R$ and $u:(-\eps,+\eps)\to Y_*$, which we denote by $\lm(s), u(s)$, which are analytic as functions of $s$,
and satisfy, $\lm(0)=\lm_*$, $u(0)=0$ and in particular
\beq\label{210216.1}
\psi(s):=\psi_*+s\phi^*_1+ u(s),\quad \mbox{satisfies}\quad \psi(0)=\psi_*,
\eeq
and is a solution of
\beq\label{3.1}
\graf{-\Delta \psi(s) =\rh_{\sscp \lm(s)}(\psi(s)) \quad \om \\ \psi(s) = 0 \qquad \qquad\pa\om}
\eeq

for any $s\in(-\eps,\eps)$. Moreover,
$$
u^{'}(0)=0,\quad \psi^{'}(0)=\phi^*_1,\quad \lm^{'}(0)=0.
$$

For later use we define

\beq\label{crux}
a_*:= {\|\psi_*\|_{\lm_*}}\al_{\lm_*,\sse}(\phi^*_1)\equiv
 <[\psi_*]_{\sscp \lm_*,0},\phi^*_1>_{\sscp \lm_*}\neq 0,
\eeq

which do not vanish because of Proposition \ref{prop1} $(i)$.

\bigskip

To understand the qualitative behavior of $\lm(s)$, we start with the following,
\bpr\label{prcrux} Let $\sg_{1,s}:=\sg_1(\lm(s),\psi(s))$, $s\in (-\eps,\eps)$ be the first eigenvalue of \eqref{2.1}.
Then, by choosing a smaller $\eps>0$ if necessary, $\lm^{'}(s)\neq 0$  in $(-\eps,\eps)\setminus\{0\}$ and $\sg_{1,s}$ and $\lm^{'}(s)$ have
the same sign in $(-\eps,\eps)\setminus\{0\}$. In particular $a_*>0$.
\epr
\proof
We first observe that $\psi^{'}(s)$ satisfies
\beq\label{2505.1}
\graf{-\Delta \psi^{'}(s) =\lm(s)\rh_{\sscp \lm(s)}[\psi^{'}(s)]_{\sscp \lm(s),0}+\lm^{'}(s)\rh_{\sscp \lm(s)}[\psi(s)]_{\sscp \lm(s),0}\quad \om \\ \\
\psi^{'}(s) = 0 \qquad \qquad\pa\om}
\eeq
and that, since $\lm$ is analytic, and by choosing a smaller $\eps$ if necessary, then
$\lm^{'}(s)\neq 0$ for $s\in (-\eps,\eps)\setminus \{0\}$. In particular we can choose $\eps$ so small to guarantee that the first eigenvalue $\sg_{1,s}$
is simple, see \cite{but}. Letting $\phi_{s}$ be the eigenfunction of the first eigenvalue $\sg_{1,s}$, we can
use \rife{lineq}, \rife{2505.1} and integrate by parts to conclude that,
\beq\label{lamq}
\sg_{1,s}\ino \rh_{\sscp \lm(s)}[\phi_s]_{\sscp \lm(s),0}\psi^{'}(s) =
\lm^{'}(s)\ino\rh_{\sscp \lm(s)}[\phi_s]_{\sscp \lm(s),0}\psi(s),\quad s\in (-\eps,\eps)\setminus \{0\}.
\eeq

As a consequence of \rife{lamq}, we obtain the following asymptotic relation,

\beq\label{crux4}
\frac{\sg_{1,s}}{\lm^{'}(s)}=
\frac{\ino\rh_{\sscp \lm(s)}[\phi_s]_{\sscp \lm(s),0}\psi(s)}{\ino \rh_{\sscp \lm(s)}[\phi_s]_{\sscp \lm(s),0}\psi^{'}(s)}=
\frac{\ino\rh_{*}[\phi^*_1]_{\sscp \lm_*,0}\psi_*}{\ino \rh_{*}[\phi^*_1]_{\sscp \lm_*,0}\phi^*_1}\to
\frac{\ino\rh_{*}[\phi^*_1]_{\sscp \lm_*,0}\psi_*}{\left\|\phi^*_1\right\|_{\sscp \lm_*}^2},\;\mbox{ as }\;s\to 0^{\pm}.
\eeq

Since $\lm$ is monotonic increasing along $\mathcal{G}_{\lm_*}$, we can assume without loss of generality
that $(\lm(s), \psi(s))\in\mathcal{G}_{\lm_*}$ in $(-\eps,0)$ and $\lm^{'}(s)>0$ in $(-\eps,0)$. In this situation, by construction,
we also have $\sg_{1,s}>0$ in $(-\eps,0)$,
which, in view of Proposition \ref{prop1} $(i)$ and \rife{crux4}, implies that

\beq\label{crux5}
a_*=\ino\rh_{*}[\psi_*]_{\sscp \lm_*,0}\phi^*_1 \equiv \ino\rh_{*}[\phi^*_1]_{\sscp \lm_*,0}\psi_*>0.
\eeq
Therefore, for $\eps>0$ small enough, $\sg_{1,s}$ and $\lm^{'}(s)$ have the same sign in $(-\eps,\eps)\setminus\{0\}$, as claimed.
\finedim

\brm\label{rem64} {\it The relation \eqref{lamq} and its consequences will play a crucial role in the following.
It is worth to point out that a fundamental role in the derivation of \eqref{lamq} is due to the analyticity of $\lm$. If we miss this property, then
we could pursue
other well known arguments in bifurcation analysis, see \cite{Amb}, and for example conclude that,
\beq\label{crux3}
\lm^{'}(s)=-\frac{b_*}{a_*}s+\bop(s),\quad s\to 0^{\pm},
\eeq
where,
$$
b_*=
(\lm_*)^2\left(<(\phi^*_1)^3>_{\sscp \lm_*}-<(\phi^*_1)^2>_{\sscp \lm_*}<\phi^*_1>_{\sscp \lm_*}\right).
$$
However we could not exclude that $b_*=0$. This case can be ruled out for some specific nonlinearities, but it is hard to fix in general,
see \cite{Amb} for further details.  If $\lm$ were not analytic, then the same problem could arise in principle for higher order terms in
the expansion of $\lm^{'}$, which (in principle) could even vanish in a full right neighbourhood of $s=0$. This would be a major problem
for our purposes, since it would affect the strict monotonicity of the energy along the branch.
}
\erm

\bigskip
\bigskip
By taking a smaller $\eps>0$ if necessary, we may define $\lm_\eps=\lm(\eps)$, and
$$
E^{(-)}(\lm):= \mathcal{E}(\rl(\pl)),\quad (\lm,\pl)=(\lm(s),\psi_{\ssb(s)}),\;s\in[0,\eps),
$$
and observe that, setting $\vb=\frac{d\pl}{d\lm}$,
$$
<[\psi(s)]_{\ssb(s),0} v_{\ssb(s)}>=\frac{1}{\lm^{'}(s)}<[\psi(s)]_{\ssb(s),0} \psi^{'}(s)>=
\frac{1}{\lm^{'}(s)}<[\psi_*]_{\sscp \lm_*,0} \phi_1^{*}>(1+\bop(1)),
$$
as $s\to 0^{\pm}$. Therefore, in view of Proposition \ref{prcrux}, we conclude that, either,\\
${\bf (I)_{\lm_*}}$ $\lm^{'}(s)>0$ and $\sg_{1,s}>0$ in $(0,\eps)$: therefore, in this case, $\lm_\eps>\lm_*$, $\lm^{'}(s)\searrow 0^+$ as $s\to 0^+$ and,
by Proposition \ref{pr2.1} $(ii)$, we find $\frac{d}{d\lm}E^{(-)}(\lm)>0$ in $(\lm_*,\lm_\eps)$,  and, in view of \eqref{8.12.11},
$$
\frac{d}{d\lm}E^{(\pm)}(\lm)=<\psi^2_{\ssb(s),0}>+\lm<[\psi(s)]_{\ssb(s),0} v_{\ssb(s)}>=
$$
$$
<[\psi_*]^2_{\ssb_*,0}>+\frac{\lm_*}{\lm^{'}(s)}<[\psi_*]_{\sscp \lm_*,0} \phi_1^{*}>(1+\bop(1))\to +\ii,\;\mbox{as }s\to 0^{\mp},
$$
or else,\\
${\bf (D)_{\lm_*}}$ $\lm^{'}(s)<0$ and $\sg_{1,s}<0$ in $(0,\eps)$. By choosing a smaller $\eps$ if necessary
we can assume that $0$ is not an eigenvalue of \eqref{2.1}: therefore, in this case, $\lm_\eps<\lm_*$, $\lm^{'}(s)\nearrow 0^-$  as $s\to 0^+$, and then, in view of \eqref{8.12.11},
$$
\frac{d}{d\lm}E^{(\pm)}(\lm)=<\psi^2_{\ssb(s),0}>+\lm<[\psi(s)]_{\ssb(s),0} v_{\ssb(s)}>=
$$
$$
<[\psi_*]^2_{\ssb_*,0}>+\frac{\lm_*}{\lm^{'}(s)}<[\psi_*]_{\sscp \lm_*,0} \phi_1^{*}>(1+\bop(1))\to \pm \ii,\;\mbox{as }s\to 0^{\mp},
$$
that is, by taking a smaller $\eps$ if necessary, $\frac{d}{d\lm}E^{(-)}(\lm)<0$ for $\lm\in (\lm_\eps,\lm_*)$.

\bigskip
\bigskip
In other words, in a right neighbourhood of $E_*$, say $[E_*,E_\eps)$ is well defined the inverse of $E^{(-)}(\lm)$, say $\lm^{(-)}(E)$, which is
continuous and monotone in $[E_*,E_\eps)$. In particular the function $\lm_{\mathcal{E}}:[0,E_\eps)\to \R$, defined by
$\lm_{\mathcal{E}}=\lm^{(+)}$ in $[0,E_*)$, $\lm_{\mathcal{E}}=\lm^{(-)}$ in $[E_*,E_\eps)$ is continuous and differentiable in
$(0,E_\eps)$ and, if ${\bf (I)_{\lm_*}}$ holds, it satisfies,
$$
\lim\limits_{E\nearrow E_*}\frac{d\lm^{(+)}(E)}{dE}= 0^+=\frac{d\lm_{\mathcal{E}}(E_*)}{dE}= 0^+=
\lim\limits_{E\searrow E_*}\frac{d\lm^{(-)}(E)}{dE},
$$

while if ${\bf (D)_{\lm_*}}$ holds, it satisfies,

$$
\lim\limits_{E\nearrow E_*}\frac{d\lm^{(+)}(E)}{dE}= 0^+=\frac{d\lm_{\mathcal{E}}(E_*)}{dE}= 0^-=
\lim\limits_{E\searrow E_*}\frac{d\lm^{(-)}(E)}{dE}.
$$

So, either ${\bf (I)_{\lm_*}}$ holds, and then $\lm_{\mathcal{E}}$ has an increasing flex with horizontal tangent at $E=E_*$, or

${\bf (D)_{\lm_*}}$ holds, and then $\lm_{\mathcal{E}}$ has a local maximum at $E=E_*$.\\

Therefore we have two possibilities.\\
If ${\bf (D)_{\lm_*}}$ holds, then $\mathcal{G}_{\lm_*}$ can be continued in a small enough left neighbourhood of
$\lm_*$ to an analytic branch which we denote by $\mathcal{G}^*_{\lm_{\sscp D}}$, with $\lm_{\sscp D}<\lm_*$, such that
$\lm_{\mathcal{E}}$ is strictly decreasing in $(E_*,E_{\lm_{\sscp D}})$,
where $E_{\lm_{\sscp D}}=\mathcal{E}(\rh_{\sscp \lm_{\sscp D}}(\psi_{\sscp \lm_{\sscp D}}))$. In other words, in this case
$\lm_{\mathcal{E}}$ is strictly decreasing as the
energy increases along $\mathcal{G}^*_{\lm_{\sscp D}}$. Moreover, 
since $0$ is not an eigenvalue of \eqref{2.1}, then we can use Lemma \ref{analytic}-$(ii)$, and so in particular Proposition \ref{enrgm} can be applied to conclude that
$\mathcal{G}^*_{\lm_{\sscp D}}\subset \mathbb{S}_\om$ and in particular that, $\mathcal{G}^*_{\lm_{\sscp D}}\setminus \ov{ \mathcal{G}_{\lm_*}}=(\lm_{\mathcal{E}}(E),\psi_{\lm_{\mathcal{E}}(E)}),E\in (E_*,E_{\lm_{\sscp D}})$,
for some
$E_{\lm_{\sscp D}}>E_*$, with $\frac{d\lm_{\mathcal{E}}(E)}{d E}<0$,$\fo\,E\in (E_*,E_{\lm_{\sscp D}})$.
Then we define $E_d=E_{\lm_{\sscp D}}$ and we have $E_m=E_*$,
which, together with Remark \ref{formulas}, readily yields the claims of Theorem \ref{thm2} $(i)-(ii)-(iii)$ in the case where \rife{convex0} holds.\\

Otherwise, if ${\bf (I)_{\lm_*}}$ holds, then $\mathcal{G}_{\lm_*}$ can be continued in a right neighbourhood of
$\lm_*$ to an analytic branch which we denote by $\mathcal{G}^*_{\mu}$, with $\mu>\lm_*$, such that
$\lm_{\mathcal{E}}$ is strictly increasing. Then $\lm_{\mathcal{E}}$ has an increasing flex with horizontal tangent at $E=E_*$ and
in this case we define,
$$
\lm^{(1)}_I:=\sup\{\mu>\lm_*\,:\, \sg_1(\lm,\pl)>0,\fo\; (\lm,\pl)\in \mathcal{G}^*_{\mu}\},
$$
and we see from Proposition \ref{pr2.1} that $\mathcal{E}$ is strictly increasing in $(\lm_*,\lm^{(1)}_I)$. Therefore
Proposition \ref{enrgm} shows that $\mathcal{G}^*_{\lm^{(1)}_I}\subset\mathbb{S}_\om$ and
that the analytic branch
${\mathcal{G}^*_{\lm^{(1)}_I} }$, is indeed a branch of entropy maximizers which can be parametrized by the
energy, that is
$$
\mathcal{G}^*_{\lm^{(1)}_I} = \left\{ (\lm_{\mathcal{E}}(E),\psi_{\lm_{\mathcal{E}}(E)}),E\in (0,E^{(1)}_I) \right\},
$$
where $E^{(1)}_I=\lim\limits_{\lm\nearrow\lm^{(1)}_I}\mathcal{E}(\rl(\pl))$. Observe that any entropy maximizer is non degenerate by assumption.
In particular, if ${\bf (I)_{\lm_*}}$ holds, and since $\om$ is strictly starshaped, then \rife{poho} and Lemma \ref{ebound}
show that $\lm^{(1)}_I<+\ii$, $E^{(1)}_I<+\ii$ and $\pl\to\psi_{\lm^{(1)}_I}$ with $(\lm^{(1)}_I,\psi_{\lm^{(1)}_I})\in \mathbb{S}_\om$. At this point,
the local analysis worked out in Proposition \ref{prcrux} at $\lm=\lm_*$, works as it stands at $\lm=\lm^{(1)}_I$ as well, and we
can continue  the branch $\mathcal{G}^*_{\lm^{(1)}_I}$ as above, either in a right
neighbourhood of $\lm^{(1)}_I$ (a case which we denote by ${\bf (I)^{(2)}_{\lm^{(1)}_I}}$) or in a left neighbourhood
of $\lm^{(1)}_I$ (a case which we denote by ${\bf (I.D)_{\lm^{(1)}_I}}$).
Clearly, if ${\bf (I.D)_{\lm^{(1)}_I}}$, occurs then we may conclude the proof as in case ${\bf (D)_{\lm_*}}$ above by setting $E_1=E_*$,
$E_m=E^{(1)}_I$ with a suitable  $E_d>E_m$. On the other side,
if ${\bf (I)^{(2)}_{\lm^{(1)}_I}}$ occurs, then $\lm_{\mathcal{E}}$ has an increasing flex with horizontal tangent at $E=E^{(1)}_I$ and by
Proposition \ref{enrgm} we can continue the branch to a branch of entropy maximizers parametrized by the energy
in a right neighbourhood of $\lm^{(1)}_I$. Therefore, if  ${\bf (I)^{(2)}_{\lm^{(1)}_I}}$ occurs, then
we can iterate the above argument and assume without loss of generality that after $n$ steps we have found $E^{(1)}_I,\cdots,E^{(n)}_I$
which are all flex with horizontal tangent corresponding to the values $\lm^{(1)}_I,\cdots,\lm^{(n)}_I$ such that
${\bf (I)^{(k)}_{\lm^{(k-1)}_I}}$ holds for $k=1,\cdots,n$, and that $\frac{d\lm_{\mathcal{E}}}{dE}>0$ if $E\neq E^{(k)}_I$, so that
in particular ${\bf (I^{(k)}.D)_{\lm^{(k-1)}_I}}$ never occurs for $k=1,\cdots,n$.\\
We claim that there exists $n_0\geq 1$ such that  ${\bf (I^{(n_0)}.D)_{\lm^{(n_0-1)}_I}}$ holds. In fact, arguing by contradiction,
if this was not the case, then
there should exist sequences $E^{(n)}_I$ and $\lm^{(n)}_I$, $n\in \N$, such that ${\bf (I)^{(n)}_{\lm^{(n-1)}_I}}$ would hold for any
$n$ and clearly ${\bf (I^{(n)}.D)_{\lm^{(n-1)}_I}}$ never occurs. Since $\om$ is strictly starshaped, then \rife{poho}
show that there exists $\widehat{\lm}>8\pi$
such that if $(\lm,\pl)$ solves $\prl$, then $\lm\leq \widehat{\lm}$. Then, by arguing as in Lemma \ref{ebound}, we see that
$\sup\limits_{n\in\N}E^{(n)}_I\leq \ov{E}<+\ii$ as well.
Thus, passing to a subsequence if necessary, we could find $\widehat{\lm}\geq \lm_I=\lim\limits_{n} \lm^{(n)}_I$,
$E_I=\lim\limits_{n} E^{(n)}_I$, and $\psi_I$
which solves $\prl$ with $\lm=\lm_I$ such
that $\psi_{\sscp \lm^{(n)}_I}\to \psi_I$ in $C^{2,\al}_0(\om)$-norm. Therefore $(\lm_I,\psi_I)\in \mathbb{S}_\om$, and
since $\sg_{1}(\lm^{(n)}_I, \psi_{\sscp  \lm^{(n)}_I})=0$, then we also have $\sg_{1}(\lm_I, \psi_I)=0$. In other words we could
apply once more the analytic implicit function theorem, and obtain an analytic extension of $\lm_{\mathcal{E}}$ in a neighbourhood of
$E_I$. However this is impossible since then we should have $\frac{d\lm_{\mathcal{E}}}{dE}\neq 0$ in a small enough
neighbourhood of $E_I$, while we know that $\frac{d\lm_{\mathcal{E}}(E^{(n)}_I)}{dE}=0$, $\fo\,n\in \N$. As a consequence, we have found
$n_0\geq 1$ such that  ${\bf (I^{(n_0)}.D)_{\lm^{(n_0-1)}_I}}$ holds and we conclude the proof as above by setting
$E_1=E_*$, $E_m=E^{(n_0)}_{I}$ with a suitable $E_d>E_m$. This fact together with Remark \ref{formulas}
concludes the proof of Theorem \ref{thm2} $(i)-(ii)-(iii)$ in this case as well. In particular we find,
$$
\frac{d^2S(E)}{d E^2}=-\frac{d\lm_{\mathcal{E}}(E)}{d E}\leq 0,\fo\,E\in(E_0,E_m],
$$
where the equality sign holds at $E=E_m$ and possibly at each one of the finitely many flex and obviously,
$$
\frac{d^2S(E)}{d E^2}=-\frac{d\lm_{\mathcal{E}}(E)}{d E}> 0,\fo\,E\in(E_m,E_d),
$$
as claimed.
\finedim

\bigskip
\bigskip

\section{\bf The proof of Theorem \ref{thm3}.}\label{sec7}
In this section we present the proof of Theorem \ref{thm3}.

\bigskip

{\bf The proof of Theorem \ref{thm3}.} $(i)-(ii)$\\
The starting point of our analysis is the proof of Theorem \ref{thm2} and we adopt the same notations used there.
We recall that
if ${\bf (D)_{\lm_*}}$ holds, then $\mathcal{G}_{\lm_*}$ can be continued in a small enough left neighbourhood of
$\lm_*$ to a smooth branch which we denoted by $\mathcal{G}^*_{\lm_{\sscp D}}$, with $\lm_{\sscp D}<\lm_*$, such that
$\lm_{\mathcal{E}}$ is strictly decreasing in $(E_*,E_{\lm_{\sscp D}})$, where
$E_{\lm_{\sscp D}}=\mathcal{E}(\rh_{\sscp \lm_{\sscp D}}(\psi_{\sscp \lm_{\sscp D}}))$.
In this case the energy is strictly increasing along $\mathcal{G}^*_{\lm_{\sscp D}}$ and we can apply Proposition \ref{enrgm}, to conclude that
$\mathcal{G}^*_{\lm_{\sscp D}}\subset \mathbb{S}_\om$ and in particular that $\mathcal{G}^*_{\lm_{\sscp D}}\setminus \ov{ \mathcal{G}_{\lm_*}}=(\lm_{\mathcal{E}}(E),\psi_{\lm_{\mathcal{E}}(E)}),E\in (E_*,E_{\lm_{\sscp D}})$,
for some
$E_{\lm_{\sscp D}}>E_*$, with $\frac{d\lm_{\mathcal{E}}(E)}{d E}<0$,$\fo\,E\in (E_*,E_{\lm_{\sscp D}})$.
Clearly, by taking a smaller $E_{\lm_{\sscp D}}$ if necessary,
by continuity we have $\sg_1(\lm_{\mathcal{E}}(E),\psi_{\lm_{\mathcal{E}}(E)})<0$, $\fo\,E\in (E_*,E_{\lm_{\sscp D}})$ and we define,
$$
E^{(1)}_D:=\sup\left\{E>E_*\,:\, \sg_1(\lm_{\mathcal{E}}(E),\psi_{\lm_{\mathcal{E}}(E)})<0\mbox{ and }
\frac{d\lm_{\mathcal{E}}(E)}{d E}<0\right\}.
$$

We have only two possibilities: either $E^{(1)}_D<+\ii$ or $E^{(1)}_D=+\ii$.\\

{\bf CASE 1: $E^{(1)}_D=+\ii$.}\\
In this case $\lm_{\mathcal{E}}$ is strictly decreasing while the energy increases along the branch in $(E_*,+\ii)$.
Then, by Proposition \ref{enrgm}, we can
extend $\lm_{\mathcal{E}}$ to $\lm^{(\ii)}_{\mathcal{E}}$ on $(E_0,+\ii)$ and conclude
that $\mathcal{G}_{\ii}:=\{(\lm^{(\ii)}_{\mathcal{E}}(E),\psi_{\sscp \lm^{(\ii)}_{\mathcal{E}}(E)}),\,E\in (E_0,+\ii)\}\subset \mathbb{S}_\om$ as desired.

\bigskip

{\bf CASE 2: $E^{(1)}_D<+\ii$}.\\
We are going to show that the assumption $\sg_1<0$ together with ${\bf (H1)},{\bf (H2)}$, rule out this case. Let
$$
\lm^{(1)}_D:=\lim\limits_{E\nearrow E^{(1)}_D}\lm_{\mathcal{E}}(E).
$$
Since $\mathcal{E}(\rl(\pl))\leq E^{(1)}_D$ for any $(\lm,\pl)\in\mathcal{G}^*_{\lm^{(1)}_D}$ with $\lm\in (\lm^{(1)}_D,\lm_*)$, then
we have $\|\pl\|_\ii\leq C$, for any $\lm\in [\lm^{(1)}_D,\lm_*]$ for a suitable constant $C>0$ depending only on $\lm_*$ and $E^{(1)}_D$,
see Lemma 2.1 in \cite{BdM2} for a proof
of this fact. Therefore, by a standard bootstrap argument, we can assume that $\pl\to\psi_{\lm^{(1)}_D}$, as $E\nearrow E^{(1)}_D$, where
$\psi_{\lm^{(1)}_D}$ is a solution of $\prl$ with $\lm=\lm^{(1)}_D$.\\
Since solutions of $\prl$ are unique for $\lm\leq 8\pi$, then $\lm^{(1)}_D> 8\pi$.
In fact, if $\lm^{(1)}_D\leq 8\pi$, then, since $E^{(1)}_D>E_*>E_{8\pi}$, and in view of Theorem \ref{thm1} $(i)$, we would have a
second solution of $\prl$ for any such $\lm$, which is impossible by the uniqueness result in \cite{BLin3}.\\
Clearly $\ov{\mathcal{G}(E^{(1)}_D)}\subset\mathbb{S}_\om$, whence, in view of \rife{pro.16.1}
and Proposition \ref{prop1} $(i)$, we find

$$
\left.\ino\rh_{\ssb}[\phi_1]_{\ssb,0}\pl\right|_{\lm=\lm^{(1)}_D}=
\left\|[\psi_{\lm^{(1)}_D}]_{\lm^{(1)}_D,0} \right\|_{\lm^{(1)}_D}\al_{\lm^{(1)}_D,\mathcal{E}}(\phi_1)\neq 0.
$$
On the other side, $\al_{\lm^{(1)}_D,\mathcal{E}}(\phi_1)$ is a continuous
function of $\lm$ which cannot change sign in $(\lm^{(1)}_D,\lm_*]$, because of Proposition \ref{prop1} and
since $\ov{\mathcal{G}(E^{(1)}_D)}\subset\mathbb{S}_\om$ and all
the entropy maximizers are non degenerate for $E\geq E_*$. Therefore,
since $\al_{\lm_*,\mathcal{E}}(\phi_1^*)> 0$ (see \rife{crux} and Proposition \ref{prcrux}), then we conclude that,
\beq\label{cruxD}
\left.\ino\rh_{\ssb}[\phi_1]_{\ssb,0}\pl\right|_{\lm=\lm^{(1)}_D}=\left\|[\psi_{\lm^{(1)}_D}]_{\lm^{(1)}_D,0} \right\|_{\lm^{(1)}_D}\al_{\lm^{(1)}_D,\mathcal{E}}(\phi_1)>0.
\eeq
At this point, we will obtain a contradiction to the assumption $\sg_1(\lm,\pl)<0$ by showing that,
\ble\label{lem7.1}
It holds
$$
\sg_{1,1}:=\sg_{1}(\lm^{(1)}_D,\psi_{\lm^{(1)}_D})=0.
$$
\ele
\proof
If $\sg_{1,1}<0$ and if
$\frac{d\lm_{\mathcal{E}}(E^{(1)}_D)}{d E}$ were well defined and $\frac{d\lm_{\mathcal{E}}(E^{(1)}_D)}{d E}<0$, then by the
implicit function theorem and by continuous dependence on the data we could continue $\lm_{\mathcal{E}}$ in a right neighbourhood of
$E^{(1)}_D$ in such a way that $\sg_{1,1}<0$ and $\frac{d\lm_{\mathcal{E}}(E)}{d E}<0$. Therefore, arguing by contradiction, if the claim were false then
we should find, either,\\
$(a)$ $\sg_{1,1}<0$ and $\frac{d\lm_{\mathcal{E}}(E^{(1)}_D)}{d E}$ is not well defined, or,\\
$(b)$ $\sg_{1,1}<0$ and $\frac{d\lm_{\mathcal{E}}(E^{(1)}_D)}{d E}=0$.\\
In both cases,
since $(\lm^{(1)}_D,\psi_{\lm^{(1)}_D})\in \mathbb{S}_\om$, then Proposition \ref{prop1} $(ii)$
implies that $0$ is not an eigenvalue of \rife{lineq}, that is, $\sg_{2,1}:=\sg_{2}(\lm^{(1)}_D,\psi_{\lm^{(1)}_D})>0$ and $v_{\sscp \lm^{(1)}_D}$ is well defined.\\

If $(a)$ holds, since
$v_{\sscp \lm^{(1)}_D}$ is well defined, there must exist a sequence $E_n\nearrow E^{(1)}_D$ such that $\lm_n:=\lm_{\mathcal{E}}(E_n)\searrow \lm^{(1)}_D$ and
$\frac{d\lm_{\mathcal{E}}(E_n)}{d E}\to -\ii$, $\ainf$. Therefore we find,

$$
<v_{\sscp \lm_n}>=\left(\frac{d\lm_{\mathcal{E}}(E_n)}{d E}\right)^{-1}\to <v_{\sscp \lm^{(1)}_D}>=0,\ainf.
$$

To prove that this is not possible we have to use ${\bf (H1)},{\bf (H2)}$.
\ble\label{sig2}
If $\sg_{1,1}<0$ and $\mbox{\rm $\rh_{\sscp \sscp \lm^{(1)}_D}$}$ is $\mu_0$-stable with $\mu_0\leq -\frac{\lm^{(1)}_D}{ \lm^{(1)}_D+\sg_{2,1}}$,
then  $<v_{\sscp \lm^{(1)}_D}>_{\sscp \lm^{(1)}_D}<0$.
\ele
\proof
Let
$$
\psi_{1,0}:=[\psi_{\sscp \lm^{(1)}_D}]_{\sscp \lm^{(1)}_D,0}=\al_1\phi_{1,0}+\psi_{1,0}^{\perp},\quad \psi_{1,0}^{\perp}
=\sum\limits_{j=2}^{+\ii}\al_j\phi_{j,0},
$$
and
$$
v_{1,0}:=[v_{\sscp \lm^{(1)}_D}]_{\sscp \lm^{(1)}_D,0}=\beta_1\phi_{1,0}+v_{1,0}^{\perp},\quad v_{1,0}^{\perp}=
\sum\limits_{j=2}^{+\ii}\beta_j\phi_{j,0},
$$
be the Fourier expansions \rife{Four} of $\psi_{1,0}=[\psi_{\sscp \lm^{(1)}_D}]_{\sscp \lm^{(1)}_D,0}$ and
$v_{1,0}=[v_{\sscp \lm^{(1)}_D}]_{\sscp \lm^{(1)}_D,0}$ with
respect to the normalized ($\|\phi_{j,0}\|_{\sscp \lm^{(1)}_D}=1$) projections $\phi_{j,0}:=[\phi_j]_{\sscp \lm^{(1)}_D,0}$ of the eigenfunctions of \rife{lineq}.
In view of \rife{cruxD} we have,
$$
\al_1 \neq 0.
$$
To avoid cumbersome notations, let us set,
$$
\rh_{1}=\rh_{\sscp \lm^{(1)}_D},\quad \|\psi_{1,0}\|_{\sscp D}=\|\psi_{1,0}\|_{\sscp \lm^{(1)}_D},\quad
\sg_{j,1}=\sg_{j}(\lm^{(1)}_D,\psi_{\sscp \lm^{(1)}_D}).
$$

By using \rife{1b1} with \rife{lineq9}, we obtain,
\beq\label{lamq3}
\sg_{j,1}\ino \rh_{1}\phi_{j,0}v_{1,0} =
\ino\rh_{1}\phi_{j,0}\psi_{1,0},\mbox{ that is }\beta_j=\frac{\al_j}{\sg_{j,1}},\quad j\in\N,
\eeq
and then we find,
\beq\label{event}
<\psi_{1,0}^{\perp},v_{1,0}^{\perp}>_{\sscp \lm^{(1)}_D}=\sum\limits_{j=2}^{+\ii} \al_j\beta_j=
\sum\limits_{j=2}^{+\ii} \sg_{j,1}(\beta_j)^2\geq \sg_{2,1}<(v_{1,0}^{\perp})^2>_{\sscp \lm^{(1)}_D}\geq 0.
\eeq

At this point we argue by contradiction and suppose that,
$$
<v_{\sscp \lm^{(1)}_D}>_{\sscp \lm^{(1)}_D}\geq 0,
$$
so that, by \rife{2b1} and \rife{event}, we find,

$$
<v_{\sscp \lm^{(1)}_D}>_{\sscp \lm^{(1)}_D}=\al_1^2+\lm^{(1)}_D\al_1\beta_1+y_1^2\geq 0,
\mbox{ where }y_1^2=\|\psi_{1,0}^{\perp}\|^2_{\sscp \lm^{(1)}_D}+ \lm^{(1)}_D
<\psi_{1,0}^{\perp},v_{1,0}^{\perp}>_{\sscp \lm^{(1)}_D}>0,
$$
and so we conclude that,
\beq\label{zero}
\al_1^2+ \frac{\al_1^2\lm^{(1)}_D}{\sg_{1,1}}+y_1^2\geq 0,\mbox{ that is }\sg_{1,1}\leq  -\frac{\lm^{(1)}_D \al_1^2}{\al_1^2+y_1^2}.
\eeq

On the other side, for $\lm=\lm^{(1)}_D$, we can test \rife{080216.1nnn} with $\varphi_{\sscp T}=\phi_{1,0}-\al_{1,\mathcal{E}}\psi_{1,0}$ where,
$$
\al_{1,\mathcal{E}}=
\frac{1}{\|\psi_{1,0}\|_{\sscp D}^2}\ino \rh_{1}\psi_{1,0}\phi_{1,0}=\frac{\al_1}{\|\psi_{1,0}\|_{\sscp D}^2},
$$
and then, by using \rife{lineq9} and \rife{lineqe}, after a straightforward calculation obtain that,
\beq\label{crux21}
\frac{\al_1^2}{\|\psi_{1,0}\|^4_{\sscp D}}\mathcal{A}_{\lm^{(1)}_D}(\psi_{1,0})+
\frac{2\sg_{1,1}}{\lm^{(1)}_D+\sg_{1,1}}\frac{\al_1^2}{\|\psi_{1,0}\|^2_{\sscp D}}-
\frac{\sg_{1,1}}{\lm^{(1)}_D+\sg_{1,1}}\leq \mu_0\left(1-\frac{\al_1^2}{\|\psi_{1,0}\|^2_{\sscp D}}\right).
\eeq

We first solve \rife{crux21} with respect to $\sg_{1,1}$, and, recalling that $\lm^{(1)}_D+\sg_{1,1}>0$, conclude that,
\beq\label{crux22}
\left(\|\psi_{1,0}\|^4_{\sscp D}-2\|\psi_{1,0}\|^2_{\sscp D}\al_1^2-\mathcal{A}_{\lm^{(1)}_D}(\psi_{1,0})\al_1^2\right)\sg_{1,1}\geq
\eeq
$$
\lm^{(1)}_D\mathcal{A}_{\lm^{(1)}_D}(\psi_{1,0})\al_1^2-\sg_0(\|\psi_{1,0}\|^4_{\sscp D}-\|\psi_{1,0}\|^2_{\sscp D}\al_1^2),
$$

where
\beq\label{crut}
\sg_0:=(\lm^{(1)}_D+\sg_{1,1})\mu_0\leq -\lm^{(1)}_D \frac{\lm^{(1)}_D+\sg_{1,1}}{\lm^{(1)}_D+\sg_{2,1}}.
\eeq

Clearly, by using \rife{lineq9} and \rife{lineqe}, we see that,
$$
\mathcal{A}_{\lm^{(1)}_D}(\psi_{1,0})=\al_1^2\frac{\lm^{(1)}_D}{\lm^{(1)}_D+\sg_{1,1}}+\mu^2_2-\|\psi_{1,0}\|^2_{\sscp \lm^{(1)}_D}
\equiv -\al_1^2\frac{\sg_{1,1}}{\lm^{(1)}_D+\sg_{1,1}}-\nu^2_2,
$$
where,
$$
\mu^2_2=\sum\limits_{j=2}^{+\ii}\al_j^2\frac{\lm^{(1)}_D}{\lm^{(1)}_D+\sg_{j,1}},
\quad \nu^2_2=\sum\limits_{j=2}^{+\ii}\al_j^2\frac{\sg_{j,1}}{\lm^{(1)}_D+\sg_{j,1}},
$$
and thus \rife{crux22} is in turn equivalent to,
$$
\left(\|\psi_{1,0}\|^4_{\sscp D}-\al_1^2\|\psi_{1,0}\|^2_{\sscp D}-\al_1^2\mu_2^2\right)\sg_{1,1}\geq
-\lm^{(1)}_D\al_1^2\nu_2^2-\sg_0(\|\psi_{1,0}\|^4_{\sscp D}-\|\psi_{1,0}\|^2_{\sscp D}\al_1^2).
$$

It is easy to check that $\|\psi_{1,0}\|^4_{\sscp \lm^{(1)}_D}-\al_1^2\|\psi_{1,0}\|^2_{\sscp \lm^{(1)}_D}-\al_1^2\mu_2^2>0$, whence we find,
$$
\sg_{1,1}\geq \frac{-\lm^{(1)}_D\al_1^2\nu_2^2-\sg_0(\|\psi_{1,0}\|^4_{\sscp D}-\|\psi_{1,0}\|^2_{\sscp D}\al_1^2)}
{\|\psi_{1,0}\|^4_{\sscp D}-\al_1^2\|\psi_{1,0}\|^2_{\sscp D}-\al_1^2\mu_2^2},
$$

which we can use together with \rife{zero} to conclude that,
\beq\label{crux45}
\frac{-\lm^{(1)}_D\al_1^2\nu_2^2-\sg_0(\|\psi_{1,0}\|^4_{\sscp D}-\|\psi_{1,0}\|^2_{\sscp D}\al_1^2)}
{\|\psi_{1,0}\|^4_{\sscp D}-\al_1^2\|\psi_{1,0}\|^2_{\sscp D}-\al_1^2\mu_2^2}\leq -\frac{\lm^{(1)}_D \al_1^2}{\al_1^2+x_1^2}.
\eeq

At this point we are left to show that \rife{crux45} yields a contradiction. In fact, \rife{crux45} is equivalent to,
$$
\|\psi_{1,0}\|^4_{\sscp D}-\al_1^2\|\psi_{1,0}\|^2_{\sscp D}-\al_1^2\mu_2^2\leq
$$
$$
\left(\nu_2^2+\frac{\sg_0}{\lm^{(1)}_D\al_1^2}(\|\psi_{1,0}\|^4_{\sscp D}-\|\psi_{1,0}\|^2_{\sscp D}\al_1^2)\right)
\left( \al_1^2+\|\psi_{1,0}^{\perp}\|^2_{\sscp D}+ \lm^{(1)}_D <\psi_{1,0}^{\perp},v_{1,0}^{\perp}>_{\lm^{(1)}_D}\right),
$$
and since,
$$
\|\psi_{1,0}\|^4_{\sscp D}-\al_1^2\|\psi_{1,0}\|^2_{\sscp D}-\al_1^2\mu_2^2=
\al_1^2\nu_2^2+\|\psi^{\perp}_{1,0}\|^4_{\sscp D},
$$
and,
$$
<\psi_{1,0},v_{1,0}>_{\lm^{(1)}_D}=\sum\limits_{j=2}^{+\ii}\al_j{\beta_j}=\sum\limits_{j=2}^{+\ii}\frac{\al_j^2}{\sg_{j,1}},
$$
then, after a straightforward calculation, we find that,
\beq\label{crux25}
\al_1^2\nu_2^2+\|\psi^{\perp}_{1,0}\|^4_{\sscp D}\leq \al_1^2\nu_2^2+\nu_2^2\xi_2^2+
\frac{\sg_0}{\lm^{(1)}_D\al_1^2}\|\psi^\perp_{1,0}\|^2_{\sscp D}(\al_1^2+\|\psi^\perp_{1,0}\|^2_{\sscp D})
\left( \al_1^2+\xi_2^2\right),
\eeq
where,
$$
\xi_2^2=\sum\limits_{j=2}^{+\ii}\al_j^2\frac{\lm^{(1)}_D+\sg_{j,1}}{\sg_{j,1}}.
$$

By using the fact that $\|\psi^{\perp}_{1,0}\|^2_{\sscp D}=\sum\limits_{j=2}^{+\ii}\al_j^2$, we can write \rife{crux25} as follows,
$$
\al_1^2\lm^{(1)}_D
\sum\limits_{i\neq j; i,j=2}^{+\ii}\al_i^2\al_j^2\left(1-\frac{\lm^{(1)}_D+\sg_{j,1}}{\sg_{j,1}}\frac{\sg_{i,1}}{\lm^{(1)}_D+\sg_{i,1}}\right)\leq
$$
$$
\sg_0\sum\limits_{j=2}^{+\ii}\al_j^2\left(\al_1^2+\sum\limits_{j=2}^{+\ii}\al_j^2\right)\left(\al_1^2+
\sum\limits_{j=2}^{+\ii}\al_j^2\frac{\lm^{(1)}_D+\sg_{j,1}}{\sg_{j,1}}\right),
$$
which in turn, after a lengthy evaluation, is seen to be equivalent to,

\beq\label{crux99}
\al_1^2\lm^{(1)}_D
\sum\limits_{i\neq j;i,j=2}^{+\ii}\al_i^2\al_j^2\left(
1-\frac{\lm^{(1)}_D+\sg_{j,1}}{\sg_{j,1}}\frac{\sg_{i,1}}{\lm^{(1)}_D+\sg_{i,1}}+
\frac{|\sg_0|}{\lm^{(1)}_D}\frac{\lm^{(1)}_D+\sg_{j,1}}{\sg_{j,1}}+\frac{|\sg_0|}{\lm^{(1)}_D}
\right)
\leq
\eeq
$$
\sg_0\al_1^2\left(\sum\limits_{j=2}^{+\ii}\al_j^4\frac{\lm^{(1)}_D+\sg_{j,1}}{\sg_{j,1}}+\sum\limits_{j=2}^{+\ii}\al_j^4+
\al_1^2 \|\psi^{\perp}_{1,0}\|^2_{\sscp D}+ \|\psi^{\perp}_{1,0}\|^4_{\sscp D}\xi_2^2 \right)<0,
$$
where the last inequality follows from \rife{crut}. We will obtain a contradiction by showing that the sum in \rife{crux99} is positive.
In fact, let us observe that,
$$
\sum\limits_{i\neq j;i,j=2}^{+\ii}\al_i^2\al_j^2\left(
1-\frac{\lm^{(1)}_D+\sg_{j,1}}{\sg_{j,1}}\frac{\sg_{i,1}}{\lm^{(1)}_D+\sg_{i,1}}+
\frac{|\sg_0|}{\lm^{(1)}_D}\frac{\lm^{(1)}_D+\sg_{j,1}}{\sg_{j,1}}+\frac{|\sg_0|}{\lm^{(1)}_D}
\right)=
$$
$$
\sum\limits_{i\neq j;i,j=2}^{+\ii}\al_i^2\al_j^2\left(
\lm^{(1)}_D\frac{\sg_{j,1}-\sg_{i,1}}{(\lm^{(1)}_D+\sg_{i,1})\sg_{j,1}}+
\frac{|\sg_0|}{\lm^{(1)}_D}\frac{\lm^{(1)}_D+2\sg_{j,1}}{\sg_{j,1}}
\right)=
$$

$$
\sum\limits_{i< j;i,j=2}^{+\ii}\frac{\al_i^2\al_j^2}{\sg_{i,1}\sg_{j,1}} \left(
-\lm^{(1)}_D\frac{(\sg_{j,1}-\sg_{i,1})^2}{(\lm^{(1)}_D+\sg_{i,1})(\lm^{(1)}_D+\sg_{j,1})}+
|\sg_0|\frac{\lm^{(1)}_D(\sg_{j,1}+\sg_{i,1})+4\sg_{j,1}\sg_{i,1}}{\lm^{(1)}_D}
\right).
$$

Therefore, in view of \rife{crut}, and since $\sg_j=\sg_i+x_{i,j}$ for some $x_{i,j}\geq 0$ whenever $i<j$, we conclude that
the sum in \rife{crux99} takes the form,
$$
\sum\limits_{i< j;i,j=2}^{+\ii}\frac{\al_i^2\al_j^2}{\sg_{i,1}\sg_{j,1}}f_{\sg_i}(x_{i,j}),
$$
where
$$
f_{\sg}(x)=-\lm^{(1)}_D\frac{x^2}{(\lm^{(1)}_D+\sg)(\lm^{(1)}_D+\sg+x)}+
\frac{\lm^{(1)}_D+\sg_1}{\lm^{(1)}_D+\sg_2}(\lm^{(1)}_D(x+2\sg)+4\sg(\sg+x)),\;x\in [0,+\ii).
$$
By using the fact that $\sg_i \geq \sg_2$ for $i\geq 2$, elementary arguments show that $f_{\sg_i}(x)>0$ on $[0,+\ii)$ for any $i\geq 2$, whenever
$$
\sg_2\geq \frac{1}{4}\frac{\lm^{(1)}_D|\sg_1|}{\lm^{(1)}_D+\sg_1},
$$
which is the desired contradiction.
\finedim

\bigskip

In view of Lemma \ref{sig2} we see that $(a)$ cannot occur and so we are left with $(b)$.\\
However $(b)$ is neither possible, since  the identity,
$$
\left(\frac{d\lm_{\mathcal{E}}(E^{(1)}_D)}{d E}\right)^{-1}=<v_{\sscp \lm^{(1)}_D}>,
$$
holds whenever $v_{\sscp \lm^{(1)}_D}$ is well defined, that is, the unique chance to come up with a vanishing value of
$\frac{d\lm_{\mathcal{E}}(E^{(1)}_D)}{d E}$ would be $<v_{\sscp \lm^{(1)}_D}>$ to be unbounded, which is surely
impossible since we have $\sg_{1,1}<0$ and $\sg_{2,1}>0$.\\
Therefore neither $(a)$ nor $(b)$ can happen which proves Lemma \ref{lem7.1}.\finedim

\bigskip

\brm\label{stability}{\it
The nondegeneracy assumption \mbox{\rm \rife{080216.1}} is not enough to establish Lemma \ref{sig2}. Indeed,
by the same argument adopted above  \mbox{\rm (}actually it is essentially enough to take the limit $\sg_0\searrow 0$ in \mbox{\rm \rife{crux25}},
which is $\mu_0=0$ in the $\mu_0$-stability,
that is almost the same as to use the nondegeneracy \mbox{\rm )} one would obtain the condition,
$$
\nu_2^2\xi_2^2-\|\psi^{\perp}_{1,0}\|^2_{\sscp \lm^{(1)}_D}>0,
$$
which is readily seen to be always satisfied. So, no contradiction arise in this case. Therefore, it seems that
the monotonicity of the energy along a branch of nondegenerate entropy maximizers with $\sg_1<0$ may change. In other words, for the entropy to be convex in $(E_*,+\ii)$, it seems that the peak of the entropy on $\mathcal{M}_{E}$ has to be
sharp enough. The $-\frac{\lm}{\lm+\sg_2}$-stability is \mbox{\rm (}together with ${\bf (H2)}$\mbox{\rm )} a sufficient condition, and we don't know whether it is sharp or not.
However, it can be shown that a sufficient condition for an entropy maximizer \mbox{\rm $\rl$} to be $-\frac{\lm}{\lm+\sg_2}$-stable is,
$$
\lm\gamma^2_{\ssb}+\sg_1>0 \mbox{ and }\sg_2(\lm,\pl)\geq \lm\left(\frac{\lm+\gamma^2_{\ssb}\sg_1}{\lm\gamma^2_{\ssb}+\sg_1}\right),\mbox{ where }
\gamma^2_{\ssb}=\frac{\al_{\lm,1}^2}{\left\|  [\psi]_{\ssb,0}\right\|_{\ssb}^2},
$$
and $\al_{\lm,1}$ is the first Fourier coefficient \mbox{\rm \rife{Four}} of $\pl$.}
\erm

\bigskip
\bigskip

{\bf The proof of  Theorem \ref{thm3} continued.}\\
As a consequence of Lemma \ref{lem7.1} we conclude that only case $E_D^{(1)}=+\ii$ is possible and so the claims in $(i)-(ii)$ are proved.
Indeed, since $\mathcal{G}_{\ii}$ has no bifurcation points, then by the same argument used in Proposition \ref{enrgm},
we also conclude that $ \mathbb{S}_{\om} \subseteq \mathcal{G}_{\ii}$. Therefore we are just left with the proof of $(iii)$ in this situation.\\
$(iii)$ First of all, $\lm^{(\ii)}_{\mathcal{E}}$ is monotone and strictly increasing in $(0,E_*)$ and strictly decreasing in $(E_*,+\ii)$, whence $E_m=E_*$ is its unique
critical and absolute maximum point.\\
Next, letting $E_n\to +\ii$, then  $\lm^{(\ii)}_{\mathcal{E}}(E_n)$ satisfy \rife{poho} and so, arguing by contradiction, and passing to a subsequence
if necessary, we can assume without loss of generality that $\lm^{(\ii)}_{\mathcal{E}}(E_n)\to \ell>8\pi$. Indeed, we can't have
$\ell<8\pi$ since solutions of $\prl$ lying on $\mathcal{G}_{8\pi}$ are unique for $\lm\leq 8\pi$, see \cite{BLin3},
and the corresponding energies satisfy $E=E^{(+)}(\lm)\leq E_{8\pi}$ by Theorem \ref{thm1} $(i)-(ii)$.
At this point we can apply Lemma \ref{ebound} and conclude that $\ell=8\pi$, which is the desired contradiction.
In particular the uniqueness of solutions \cite{BLin3} shows that we cannot have $\lm^{(\ii)}_{\mathcal{E}}(E)=8\pi$ for some $E>E_{8\pi}$, which
implies that $\lm^{(\ii)}_{\mathcal{E}}(E)\to (8\pi)^+$, as $E\to +\ii$.\\
At this point, since obviously we can use \rife{formula}, then we immediately verify that the entropy satisfies \rife{convex}.
\finedim

\bigskip
\bigskip

\section{\bf On the global bifurcation diagram of entropy maximizers. A weaker form of Theorem \ref{thm3}.}\label{sec8}

We are concerned with the more general situation where $\sg_1=\sg_1(\lm,\pl)$ may change sign along the branch of entropy maximizers
for $E>E_d$.
We will be very sketchy about those part of the proof which is worked out as in Theorems \ref{thm2} and \ref{thm3} above.

\bte\label{thm4}
Let $\om$ be a strictly starshaped domain of second kind and suppose that
$\mathbb{S}_\om$ is pathwise connected and that
any entropy maximizer wit $E\geq E_*$ is non degenerate. Let $E_d$ be defined as in Theorem \ref{thm2}.
If, for any $(\lm,\pl)\in \{\mathbb{S}_\om \cap \mathcal{E}(\mbox{\rm $\rh$}_{\ssb}(\pl))\geq E_d\}$, it holds:\\
${\bf (H1)}$ $\mbox{\rm $\rh$}_{\ssb}$ is $\mu_0$-stable with $\mu_0\leq -\frac{\lm}{\lm+\sg_2}$ and,\\
${\bf (H2)}$ $\sg_2(\lm,\pl)\geq \frac14\frac{\lm |\sg_1|}{\lm +\sg_1}$,\\
then there exists an analytic function
$\lm^{(\ii)}_{\mathcal{E}}:(E_0,+\ii)\to (0,+\ii)$, such that:\\
$(i)$ $\left.\lm^{(\ii)}_{\mathcal{E}}\right|_{[E_0,E_d)}\equiv \lm_{\mathcal{E}}$;\\
$(ii)$ $\mathcal{G}_{\ii}:=\{(\lm^{(\ii)}_{\mathcal{E}}(E),\psi_{\sscp \lm^{(\ii)}_{\mathcal{E}}(E)}),\,E\in (E_0,+\ii)\}\subset \mathbb{S}_\om$, that is,
$\mbox{\rm $\rh$}_{\sscp \lm^{(\ii)}_{\mathcal{E}}(E)}(\psi_{\sscp \lm^{(\ii)}_{\mathcal{E}}(E)})$ is a solution
of the $\mbox{\bf (MVP)}$ at energy $E$ for each $E\in (E_0,+\ii)$;\\
$(iii)$ $\lm^{(\ii)}_{\mathcal{E}}\to (8\pi)^+$ as $E\to +\ii$ and $\lm^{(\ii)}_{\mathcal{E}}$ has possibly countably many critical points and the set of critical
points has no accumulation points;\\
$(iv)$ In particular  $\beta(E)=-\lm^{(\ii)}_{\mathcal{E}}(E)$ in $(E_0,+\ii)$ so that $S(E)$ is strictly decreasing in $(E_0,+\ii)$,
$\frac{d^2S(E)}{d E^2}=-\frac{d\lm_{\mathcal{E}}(E)}{d E}<  0$ in $I_1=(E_0, E_*)$ and there exists a set $J$, which is
either finite or at most countable but with no accumulation points, such that $E_*\in J$ and $\frac{d^2S(E)}{d E^2}=0$ in $[E_*,+\ii)$ if and only
$E\in J$.
\ete

\begin{figure}[h]
\psfrag{U}{${<\!\pl\!\!>_{\ssb}}$}
\psfrag{V}{${\lambda}$}
\psfrag{L}{$\sscp {\lambda_m}$}
\psfrag{G}{$\sscp {\lambda_*}$}
\psfrag{P}{$\sscp {8\pi}$}
\psfrag{E}{$\sscp 2E_m$}
\psfrag{F}{$\sscp 2E_*$}
\psfrag{Z}{$\sscp 2E_0$}
\includegraphics[totalheight=2in]{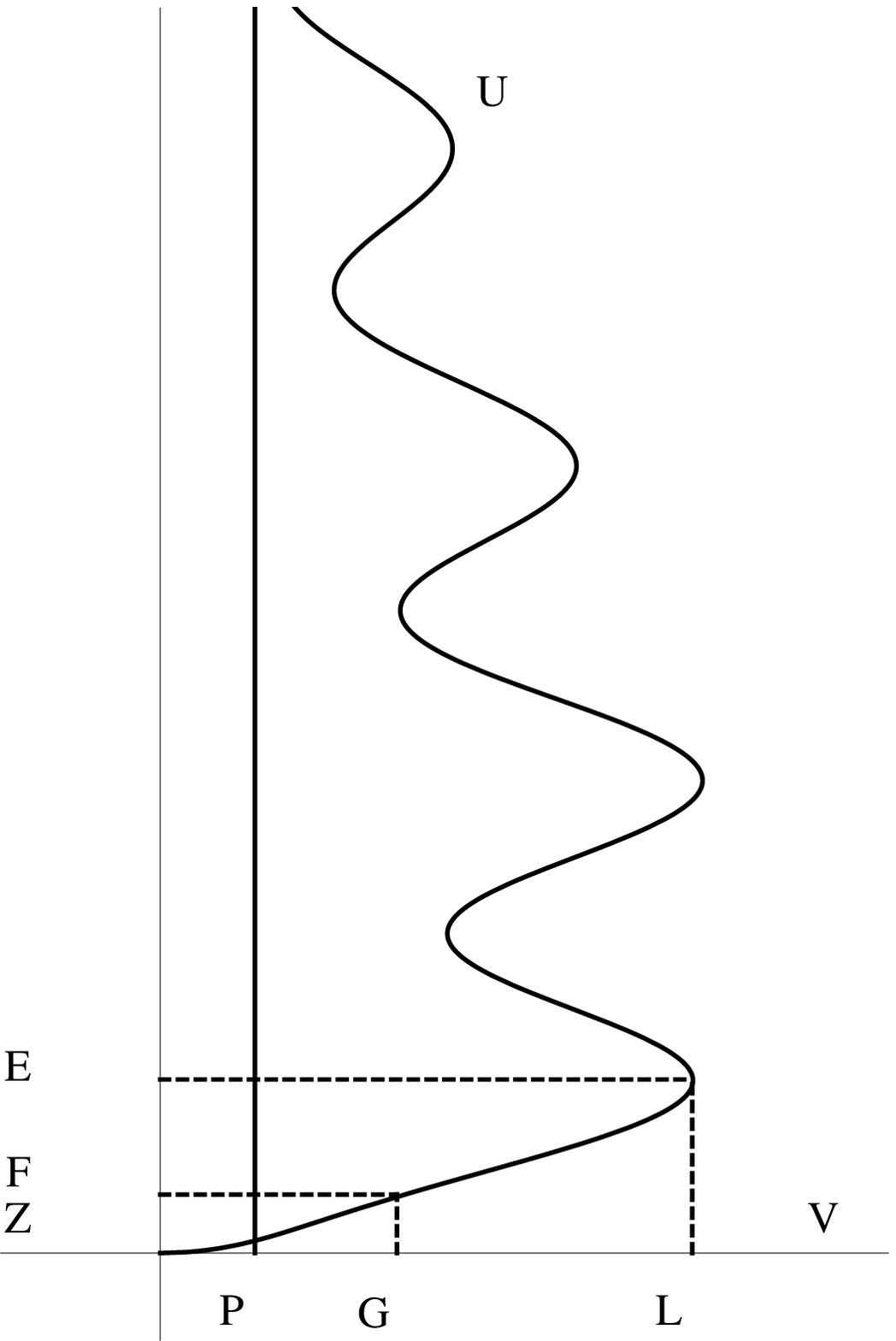}
\caption{A possible form of the global branch $\mathcal{G}_{\ii}$ in the plane $(\lm,<\!\pl\!\!>_{\ssb})$.}\label{fig4}
\end{figure}

A bifurcation diagram compatible with Theorem \ref{thm4} is depicted in fig. \ref{fig4}.

{\it The Proof of Theorem \ref{thm4}.}
Clearly Theorem \ref{thm2} holds and so we are left to prove the claims for $E\geq E_d$. In particular we can follow the proof of Theorem \ref{thm3}
and, in view of Remark \ref{rem1.11},  by defining $E^{(1)}_D$ in the same way, conclude that either CASE 1 holds that is $E^{(1)}_D=+\ii$ or CASE 2 holds, that is
$E^{(1)}_D<+\ii$. If CASE 1 holds, then the bifurcation diagram for $E\in (E_d,+\ii)$ is the same as in
Theorem \ref{thm3} and the conclusion easily follows. The situation is more subtle if CASE 2 holds. Indeed, we can define $\lm^{(1)}_D$ in the same way and conclude
via Lemma \ref{lem7.1} that,
$$
\sg_{1,1}:=\sg_{1}(\lm^{(1)}_D,\psi_{\lm^{(1)}_D})=0.
$$
The point is that now this case cannot be ruled out as above, where we had $\sg_{1,1}<0$ by assumption. Let us recall that, by \rife{cruxD}, we have
\beq\label{cruxDn}
a^{(1)}_*=\left.\ino\rh_{\ssb}[\phi_1]_{\ssb,0}\pl\right|_{\lm=\lm^{(1)}_D}=\al_{\lm^{(1)}_D,\mathcal{E}}(\phi_1)>0.
\eeq

Let $\phi^{(1)}_1$ denote the eigenfunction of $\sg_{1}(\lm^{(1)}_D,\psi_{\lm^{(1)}_D})$ and $\psi^{(1)}_1=\psi_{\lm^{(1)}_D}$.
Clearly Theorem \ref{bifurc} applies with $\widehat{\lm}=\lm^{(1)}_D$ so that, by setting,
$$
Y^{(1)}=\left\{u\in C^{2,\alpha}_0(\om)\,:\, \ino u \phi^{(1)}_1=0\right\},
$$
in a small interval $s\in(-\eps,\eps)$, there exist
$\lm^{(1)}:(-\eps,\eps) \to \R$ and $u^{(1)}:(-\eps,+\eps)\to Y^{(1)}$, which we denote by $\lm^{(1)}(s), u^{(1)}(s)$ which are analytic
functions of $s$, and satisfy, $\lm^{(1)}(0)=\lm^{(1)}_D$, $u^{(1)}(0)=0$ so that,
\beq\label{210216.1.D}
\psi(s):=\psi^{(1)}_1+s\phi^{(1)}_1+ u^{(1)}(s),\quad \mbox{satisfies}\quad \psi(0)=\psi^{(1)}_1,
\eeq
and is a solution of
\beq\label{3.1.D}
\graf{-\Delta \psi(s) =\rh_{\sscp \lm(s)}(\psi(s)) \quad \om \\ \psi(s) = 0 \qquad \qquad\pa\om}
\eeq

for any $s\in(-\eps,\eps)$. The first order derivatives of $(\lm^{(1)}(s), u^{(1)}(s))$, putting $\frac{d}{ds}:={}^{'}$, satisfy,
$$
(u^{(1)})^{'}(0)=0,\quad \psi_{\sscp \lm^{(1)}_D}^{'}(0)=\phi^{(1)}_1,\quad (\lm^{(1)})^{'}(0)=0.
$$

Moreover, in view of \rife{cruxD}, then also the argument in Proposition \ref{prcrux} applies where we just replace  $\lm_*$ with $\lm^{(1)}_D$,
and $a_*$ with $a_*^{(1)}$. In particular the analysis of  \rife{lamq} in this case is also easy, since we already now that $a_*^{(1)}>0$,
which implies that the conclusion of Proposition \ref{prcrux} hold, that is, $(\lm^{(1)})^{'}(s)\neq 0$ and has a fixed sign in $(-\eps,\eps)\setminus \{0\}$ where
$\sg_1(\lm^{(1)}(s), \psi^{(s)})$ and $(\lm^{(1)})^{'}(s)$ share the same sign.\\
At this point the discussion is very similar to that worked out in the proof of Theorem \ref{thm2} and we will be very sketchy to avoid
repetitions. Indeed, by arguing in the same way, we conclude that $\lm^{(1)}_D$ is either
a decreasing flex with horizontal tangent
or a strict local minimum (and the energy is increasing as a function of $\lm$ for $\lm>\lm^{(1)}_D$).
In particular, as in that case, the energy is always increasing, which allows  in both cases, via Proposition \ref{enrgm},
the extension of $\lm_{\mathcal{E}}$ in a right neighbourhood of $E^{(1)}_D$ as a real analytic function of $E$, in such a way that
$(\lm_{\mathcal{E}}(E),\psi_{\sscp \lm_{\mathcal{E}}(E)})\in \mathbb{S}_\om$. If $\lm^{(1)}_D$ happens to be a decreasing flex we define,
$$
E^{(1,2)}_I:=\sup\left\{E>E^{(1)}_D\,:\, \sg_1(\lm_{\mathcal{E}}(E),\psi_{\lm_{\mathcal{E}}(E)})<0,\mbox{ and }
\frac{d\lm_{\mathcal{E}}(E)}{dE}<0\right\},
$$
while if $\lm^{(1)}_D$ is a local minimum, we define,
$$
\lm^{(1)}_*:=\sup\left\{\mu>\lm^{(1)}_D\,:\,\sg_1(\lm,\pl)>0,\,\fo\,(\lm,\pl)\in\mathcal{G}_{\mu}\right\}.
$$

In the first case all the properties used in the discussion of $\lm^{(1)}_I$ are satisfied, while in the second case
all the properties used in the discussion of $\lm^*$ are satisfied so we can iterate the argument. In other words we can define
$\lm^{(\ii)}_{\mathcal{E}}$ on $(E_0,+\ii)$ as claimed in $(i)$ with the property that
$\mathcal{G}_{\lm^{(\ii)}_{\mathcal{E}}}$ is the continuation of $\mathcal{G}_{\lm_d}$ as a branch of entropy maximizers
as claimed in $(iii)$, passing possibly through a sequence of critical points of $\lm^{(\ii)}_{\mathcal{E}}$
which can be either flex or local extrema, the sign of the first eigenvalue and that of the derivative of $\lm_{\mathcal{E}}$ being always the same
along the branch. Moreover, as in the proof of Theorem \ref{thm2}, the sequence of critical points constructed in this way do not
have accumulation points.\\
The rest of the proof is the same as in Theorem \ref{thm2} and we skip it here to avoid repetitions.
\finedim

\end{document}